\documentclass[]{article}

\usepackage[utf8]{inputenc}
\usepackage{bm}
\usepackage[margin=1in]{geometry} 
\usepackage{amsmath,amsthm,amssymb,amsfonts, fancyhdr, color, comment, graphicx, environ}
\usepackage{rotating}
\usepackage{float}
\usepackage{xcolor}
\usepackage{mdframed}
\usepackage[shortlabels]{enumitem}
\usepackage{indentfirst}
\usepackage{csvsimple}
\usepackage{longtable}
\usepackage{eucal}
\usepackage{url}
\usepackage{setspace}
\usepackage{tabulary}
\usepackage{tikz}
\usepackage{amsbsy}
\usepackage{bm}
\usepackage{mathtools}
\usepackage{pythonhighlight}
\usepackage{matlab-prettifier}
\usepackage{listings}

\usepackage[flushmargin, bottom]{footmisc}
\usepackage{enumitem}
\usepackage{fancyvrb}
\usepackage{subcaption}
\usepackage{authblk}
\usepackage{letltxmacro}

\usepackage{algorithm,algpseudocode} %,algorithmic
\usepackage{hyperref}
\hypersetup{
    colorlinks=true,
    linkcolor=blue,
    filecolor=magenta,      
    urlcolor=blue,
}
\usepackage{cleveref}

\newcommand{\keywords}[1]
{
  \small	
  \textbf{Keywords:} #1.
}
\newcommand{\mathsubjclass}[2]
{
  \small	
  \textbf{Mathematical Subject Classification:} Primary: #1; Secondary: #2.
}
\newcommand{\acknowledgments}[1]
{
  \small	
  \textbf{Acknowledgments:} #1.
}

\newcolumntype{L}{>{$}l<{$}} % math-mode version of "l" column type
\newcolumntype{K}[1]{>{\centering\arraybackslash}p{#1}}

\usepackage[square,numbers]{natbib}
\title{Node Subsampling for Multilevel Meshfree Elliptic PDE Solvers}
\author[1,2]{Andrew P. Lawrence}
\author[3]{Morten E. Nielsen}
\author[2]{Bengt Fornberg}
\affil[1]{Corresponding Author, Email Address: anla5397@colorado.edu}
\affil[2]{Department of Applied Mathematics, University of Colorado Boulder, Boulder, CO 80309, USA}
\affil[3]{AAU Energy, Aalborg University, 6700 Esbjerg, Denmark}

\begin{document}
\maketitle

\begin{abstract}
    Subsampling of node sets is useful in contexts such as multilevel methods, computer graphics, and machine learning. On uniform grid-based node sets, the process of subsampling is simple. However, on node sets with high density variation, the process of coarsening a node set through node elimination is more interesting. A novel method for the subsampling of variable density node sets is presented here. Additionally, two novel node set quality measures are presented to determine the ability of a subsampling method to preserve the quality of an initial node set. The new subsampling method is demonstrated on the test problems of solving the Poisson and Laplace equations by multilevel radial basis function-generated finite differences (RBF-FD) iterations. High-order solutions with robust convergence are achieved in linear time with respect to node set size.
\end{abstract}

\\

\mathsubjclass{65N50,65N22}{65F10,65N06,65N55}

%%%%%%%%%%%%%%%%%%%%%%%%%%%%%%%%%%%%%%%%%%%%%%%%%%%%%%%%%%
\section{Introduction}

Subsampling of variable density node sets has applications in polynomial approximation, numerical integration, artificial intelligence, machine learning, multilevel methods, and computer graphics. For each of these applications, algorithms exist in 1D, 2D, and even $N$-D space, but their utility is often application specific.

%polynomical approx and quadrature
Subsampling methods have been specifically developed to choose points optimized for global polynomial approximation and numerical integration \cite{DeMarchi2015,Narayan2013,Piazzon2016,SOMMARIVA20091324}. Node sets have been optimized for global RBF collocation methods using multi-objective optimization \cite{ROQUE2014}. However, the coarse node sets these algorithms produce do not, in general, preserve the variable density of the initial, fine node sets.

%AI/ML
In the context of data driven artificial intelligence and machine learning, the process of tuning data rather than tuning model parameters has driven research on subsampling \cite{Kennard1969, meng2022, Shang2022DiversitySC, Silveira2022}. With the exception of the generalized diversity subsampling algorithm in \cite{Shang2022DiversitySC}, these algorithms are either designed for uniform subsampling or statistical learning techniques such that they are not well-suited to the preservation of variable density data sets.

%computer graphics
Research in computer graphics has led to considerable developments in the area of Poisson disk sampling which serves to subsample variable density node sets, producing resultant node sets with desirable statistical and minimum spacing properties \cite{Cook1986, Dippe1985}. The process of Poisson disk sampling is recast as a weighted sample elimination or weighted subsampling problem in \cite{Yuksel2015}. Other efforts have employed Poisson disk sampling to produce heirarchical node sets for multilevel methods using RBFs \cite{LeBorne2020}, albeit on uniform density node sets.

%Multilevel
Use of a geometric multilevel method over a variable density node set requires a subsampling routine which maintains the variable density of the original node set. Algebraic multilevel algorithms coarsen the operators themselves and the coarse levels have no intuitive geometric meaning or interpretation \cite{shapira2008matrix,trottenberg2000multigrid,falgout}. As such, coarsening methods for algebraic multilevel algorithms are not useful for geometric multilevel schemes \cite{seibold2010performance}.

%Algebraic Multigrid
Algebraic multilevel methods (AMM) provide robust and scalable linear solvers for a wide class of problems. 
They are in principle a natural choice for meshfree discretizations since the hierarchical levels are a natural byproduct of the inter-level transfer and coarse level operators. 
In the context of meshfree systems, AMM has been applied to methods that don't use RBF-FD \cite{nick2020algebraic} \cite{metsch2020algebraic} and those that do \citep{WrightJonesShankarMGMRBFFD}. For solvers which use RBFs, it has been shown that geometric multilevel methods (GMM) converge in fewer iterations \cite{WrightJonesShankarMGMRBFFD}. Additionally, the set-up time for AMM is higher overall \cite{watanabe2005comparison} \cite{volkov2016geometric}. The construction of the coarse levels themselves is higher in GMM, but that cost is reduced for a meshfree domain and motivates the need for a fast subsampling algorithm as explored in the following sections. Tests run in \cite{volkov2016geometric} demonstrate that AMM is sensitive to the mesh variation and resolution on the coarsest level. The proper choice of parameters (the strength parameter in particular) for AMM can reduce the total computation time by 15–40\%, per \cite{volkov2016geometric}. The GMMs have no such parameter sensitivity and have less sensitivity to mesh variation. According to \cite{oosterlee1998}, when using MGM and AGM as preconditioners for Krylov methods, the scheme will converge more quickly for preconditioned matrices for which the spectrum is more heavily clustered toward one. This corresponds to coefficient matrices\footnote{Those representing the application of one V-cycle of either GMM or AMM} with spectra clustered at zero. In the problems considered in \cite{WrightJonesShankarMGMRBFFD}, the spectra for MGM were more clustered around one than the compared AMM method (PyAMG \cite{Bell_PyAMG_Algebraic_Multigrid}) in all cases.
For these reasons, algebraic multigrid methods or the coarsening methods therein are not considered here.

%
%AMM are further limited to only solving problems with a single physical unknown \cite{geenen2009scalable}. The dependence of AMM on the sparsity structure of the underlying matrix tend to deteriorate when coupling vector components in higher order computations \cite{clevenger2021comparison}. While AMM are popular due to their being more easily incorporated into parallel architectures, the already growing complexity of the operators is often exacerbated by parallel implementations \cite{falgout}. While algebraic methods provide no intuitive geometric interpretations for the coarse levels, they also require no explicit knowledge of the problem geometry. Geometric methods require the user to have a sufficient knowledge of the problem geometry to place nodes appropriately, though this is frequently the case already that the user has sufficient understanding of the problem they are solving. That understanding is again important in the context of PDEs with discontinuous coefficients. Algebraic methods inherently take take these into account, while geometric methods rely on direct discretization\footnote{Such as building RBF-FD based operators at each level (as opposed to Galerkin operators) \cite{trottenberg2000multigrid}} and node placement\footnote{This emulates operator-dependent interpolation, which is the process of incorporating operator information (i.e. discontinuous coefficients) into the domain discretization \cite{falgout} \cite{trottenberg2000multigrid}} enabled by the geometric flexibility of meshfree methods to account for discontinuous coefficients.

The combination of geometric multilevel methods with meshfree solvers for partial differential equations have become increasingly popular. Meshfree methods such as radial basis function-generated finite differences (RBF-FD) discretize at scattered (quasi-uniform) nodes rather than with meshes. RBF-FD methods, in particular, allow for high geometric flexibility and can benefit from high density variation but require the underlying node sets to meet certain quality constraints in order to ensure stability and accuracy of the solution \cite{FF2015PDEwRBF,FF2015RBFPrimer,Liu2021}. Robust algorithms for generating such node sets exist \cite{FF2015NodeGen,suchde2022,vanderSande2021} and are utilised in this paper. The application of meshfree partial differential equation (PDE) solvers within a multilevel scheme requires a similarly robust algorithm for coarsening node sets \cite{Floater1996, floater1998thinning}. When implementing a multilevel algorithm, one typically starts with the initial, fine node set. Given a desired level of refinement, the task of producing a coarse node set from a fine node set can be accomplished in one of two ways: one can either select a subset of the fine node set or generate a node set that is independent of the fine node set. Many methods exist to create coarse node sets which are not subsets of the initial, fine node set \cite{drumm2008finite, suchde2019fully, suchde2021meshfree}. However, the operators to coarsen and refine between independent node sets can introduce numerical instabilities and require more memory. Alternatively, selecting a subset simplifies the coarsening and refining operators and requires less memory. The combination of a multilevel scheme with RBFs has been explored before, however, primarily on uniform (Cartesian grid) or uniformly distributed scattered node sets \cite{fasshauer2007meshfree,LeBorne2020,wendland2004scattered,Wendland2010, Wendland2017,nick2020algebraic,metsch2020algebraic}. The use of multilevel techniques on RBF-FD meshfree solvers for PDEs over variable density node sets is explored in \cite{zamolo2019novel}, however the subsampling routine used therein, based on \cite{Katz2000}, is not adjustable to coarse node sets of any size; it is limited to coarsening by factors of $1/n, n\in\mathbb{N}$. Due to this limitation, it is not considered in this paper. Though not applicable in it's original form (as it relies on information from a mesh at the fine level), an extension of the algorithm found in \cite{ha2021meshless} can be applied to meet the outlined needs for a variable density node set subsampling algorithm. However, it also suffers from inflexible coarsening factors and, as such, is not considered here. The multilevel meshfree PDE solver presented here achieves high-order solutions with robust convergence in linear time with respect to node set size.

%refining vs. coarsening
In contrast to the process of generating coarser node set from an initial fine node set, one might consider an initial coarse node set and the generation of finer node sets. Most refining algorithms use some residual function to determine if refinement should take place \cite{cavoretto2019adaptive,cavoretto2020two,zhang2017adaptive}. Most refinement techniques require user-supplied criteria in the form of a residual or principle function \cite{ling2012adaptive} at which some set of test nodes are evaluated based on the residual to determine where an how to refine. These function evaluations are an increased computational cost. Additionally, the goodness of the refinement depends heavily on the proposed test nodes. Simple ways of determining these nodes, such as using the halfway points between existing nodes \cite{cavoretto2023adaptive}, may not apply well enough to variable density node sets. On the other hand, more robust methods such as determining the Voronoi nodes \cite{behrens2001effective} or the centroids of a node and it's $K$ nearest neighbors \cite{kaennakham2019automatic} may still be ill suited for variable density refinement and introduce more significant increases in computational cost. Other refinement methods are limited to uniform refinement which is not appropriate for our current applications \cite{suchde2019fully, suchde2021meshfree}. Ultimately, refinement techniques will not be considered here.

Throughout this paper, the terms subsampling and coarsening will be used to refer to the process of selecting a subset from a collection of nodes or points. The subsampling algorithms considered in this paper are outlined in Section \ref{sec:SubAlgs}, boundary considerations for subsampling routines are covered in Section \ref{sec:BoundaryConsiderations}, numerical tests and comparisons between those presented earlier are presented in Section \ref{sec:Tests}. Additionally, Section \ref{sec:Example} includes two examples of a meshfree multilevel RBF-FD PDE solver utilizing the novel moving front node subsampling method from Section \ref{sec:MF}.

%%%%%%%%%%%%%%%%%%%%%%%%%%%%%%%%%%%%%%%%%%%%%%%%%%%%%%%%%%
\section{Subsampling Algorithms} \label{sec:SubAlgs}

This section surveys the methodology of four subsampling algorithms. In addition to a novel moving front method presented in Section \ref{sec:MF}, a weighted subsampling method based on \cite{Yuksel2015}, a method based on Poisson disk sampling, and the generalized diversity subsampling method found in \cite{Shang2022DiversitySC} are presented in Sections \ref{sec:WS}, \ref{sec:PD}, and \ref{sec:gDS} respectively.

%%%%%%%%%%%%%%%%%%%%%%%%%%%%%
\subsection{Moving Front} \label{sec:MF}

The novel subsampling algorithm presented here is a streamlined application of a 'moving front' strategy akin to those found in the node generation algorithms in \cite{vanderSande2021} and \cite{FF2015NodeGen}. Algorithm \ref{alg:MF} begins by sorting all nodes in the fine node set according to an arbitrary direction\footnote{The directional sorting and progression of the algorithm enables a cost savings in that only the nodes above the present one need to be searched}; for example, from the bottom to the top. Then, the $k$ (e.g. $k=10$) nearest neighbors to each node are determined\footnote{Achieved for a total of $N$ nodes in $O(N\log{N})$ operations by the kd-tree algorithm.}. For each node in the fine node set and working in the chosen direction, first check if the node has already been marked. If it has been marked, continue on to the next node. If it has not been marked, mark each of the $k$ nearest neighbors that is within a distance $c$ (for example $c=1.5$) of the present node's original nearest neighbor and above the present node in the sort. All marked nodes are then removed to produce the coarse node set. The moving front algorithm generalizes immediately to any number of space dimensions. A Python code for the moving front algorithm can be found in Appendix \ref{apx:MF}.

%\begin{algorithm}
%  \caption{Moving Front Algorithm} \label{alg:MF}
%  \begin{algorithmic}[1]
%  \Function{MFSub}{$X\_fine, c, K$}
%    \State $X\_fine \leftarrow \textproc{sort} (X\_fine)$
%    \State $Idx, D \leftarrow \textproc{KNearestNeighbors} (X\_fine, X\_fine, N\_neigbors = K)$
%    \For {$i = 1:\textproc{length}(X\_fine)$}
%      \If {$Idx[i,1] ~= 0$}
%        \State $ind \leftarrow \textproc{find}(D[i,2:end] < c*D[k,2])$
%        \State $ind2 \leftarrow Idx[i,ind+1]$
%        \State $ind2[ind2<i] \leftarrow [] $
%        \State $Idx[ind2,1] \leftarrow 0$
%      \EndIf
%    \EndFor
%    \State $X\_coarse \leftarrow X\_fine(Idx(:,1) ~= 0,:)$
%    \Return $X\_coarse$
%  \EndFunction
%  \end{algorithmic}
%\end{algorithm}

\begin{algorithm}
  \caption{Moving Front Algorithm} \label{alg:MF}
  \begin{algorithmic}[1]
  \Function{MFSub}{$X\_fine=\{\bm{x_1,...,x_N}\}, c, k$}
    \State Sort the nodes in $X\_fine$
    \State Find the indices and distances of the $k$ nearest neighbors for each point in $X\_fine$
    \For {$i = 1:N$}
      \If {the node $\bm{x_i}$ has not already been marked}
        \State Determine which nearest neighbors are within a radius of $c$ times the distance to the nearest neighbor
        \State Of those, determine which are 'above' $\bm{x_i}$ in the sort order
        \State Mark these nodes
      \EndIf
    \EndFor
    \State Remove the marked nodes from $X\_fine$ to produce $X\_coarse$
    \Return $X\_coarse$
  \EndFunction
  \end{algorithmic}
\end{algorithm}

It should be reiterated that intrinsic to the moving front algorithm is a directional bias. More specifically, as the algorithm proceeds across a node set, the resultant node set will differ based on the direction in which the moving front travels. The effects of this directional bias are insignificant, however, as shown in the Sections \ref{sec:subBoundSep} and \ref{sec:ExampleNumRes}.

%With respect to the computational efficiency of the algorithm, one might consider running the nearest neighbors search $N$ times, once for each point. This is because the memory usage for constructing the KD-tree and searching for nearest neighbors can be very high, and doing it once for all $N$ points can exceed the available memory of the computer. Looping over each point and finding its K nearest neighbors individually may be more memory-efficient. The memory efficiency is contrasted, however, by the overhead cost of constructing the KD-tree and searching for nearest neighbors can be significant. Performing the search once for all $N$ points will likely be more computationally efficient than repeating the process $N$ times. It is the opinion of the authors that computational efficiency is valued more.

In summary, the most efficient approach depends on the size of the input data and the available hardware resources. It's recommended to test both approaches on a subset of the data and compare their performance to determine the most efficient approach for a given problem.

%%%%%%%%%%%%%%%%%%%%%%%%%%%%%
\subsection{Weighted Subsampling} \label{sec:WS}

The weighted subsampling compared with here is based on the work presented in \cite{Yuksel2015}, but modified for variable density node sets.\footnote{A sampling example presented in \cite{Yuksel2015} should, in principle, serve to subsample variable density node sets. However, after repeated attempts, the example was not reproducible.} Each node is assigned a weight based on its distance to its nearest neighbors. The algorithm then iterates to remove the node with the highest weight, adjust the remaining weights accordingly, and repeat until the desired number of nodes remain. The code for this implementation can be found on the author's GitHub, \cite{Lawrence_GitHub}.

%%%%%%%%%%%%%%%%%%%%%%%%%%%%%
\subsection{Poisson Disk Subsampling} \label{sec:PD}

Given a radius of exclusion for each node (such that the radii of exclusion are spatially variable), a node is randomly selected from the fine node set and accepted into the coarse node set if its radius of exclusion does not overlap that of any of the previously accepted nodes. The first node in the coarse node set is chosen randomly. The radii of exclusion are the product of the distance of the nearest neighbor in the fine node set and a hyperparameter $h$, $r(x_i) = c * r_{min}(x_i)$. This algorithm is similar to the thinning method presented in \cite{suchde2022} and the Poisson thinning method presented in \cite{LeBorne2020}. The use of a nearest neighbor search to support spatially variable radii of exclusion limits the computational complexity to being no better than $O(N \log{N})$ in contrast to the $O(N)$ algorithm in \cite{bridson2007fast}. The code for this implementation can be found on the author's GitHub, \cite{Lawrence_GitHub}.

%%%%%%%%%%%%%%%%%%%%%%%%%%%%%
\subsection{Generalized Diversity Subsampling} \label{sec:gDS}

The generalized Diversity Subsampling algorithm as found in \cite{Shang2022DiversitySC} selects a subsample from the fine node set according to an arbitrary, specified distribution. The distribution utilized in this paper is a function of the distance to the nearest neighbor of each node.

%%%%%%%%%%%%%%%%%%%%%%%%%%%%%
%\subsection{Ha and Choi} \label{sec:hc}
%
%The algorithm in this section is an adaptation of the work presented in \cite{ha2021meshless}. The original authors begin with a fine meshed node set and define a set of neighbors for each point based on the mesh. Iterating through each node, if the node has not already been eliminated all of it's neighbors are eliminated. The neighbors at successive levels are defined as those nodes which had common neighbors in the level above. Given that the definition for neighboring nodes at the fine level is based on a mesh in the original method, for the present context in which no mesh exists for the fine node set, another definition must be contrived for the method to be applicable. It seems most natural to use some $k$ nearest neighbors to comprise the set of neighboring points. Doing so established $k$ as a parameter which affects not only the quality of the subsequent subsamplings, but also the 

%%%%%%%%%%%%%%%%%%%%%%%%%%%%%%%%%%%%%%%%%%%%%%%%%%%%%%%%%%
\section{Boundary Considerations} \label{sec:BoundaryConsiderations}

The purpose of this section is to illustrate the potential pitfalls that can occur if boundary nodes are included in the domain node set without being handled separately and to propose methods for overcoming those pitfalls. For the sake of demonstration, only the moving front algorithm is considered in this section. Initial node sets are generated by \cite{vanderSande2021} and nodes near the boundary have been repelled\footnote{following the repel methodology described in \cite{FF2015NodeGen}} prior to any subsampling.

%%%%%%%%%%%%%%%%%%%%%%%%%%%%%
\subsection{Subsample Boundary with Domain}
\begin{figure}[]
    \centering
    \begin{subfigure}[t]{0.45\textwidth}
        \includegraphics[width=\textwidth]{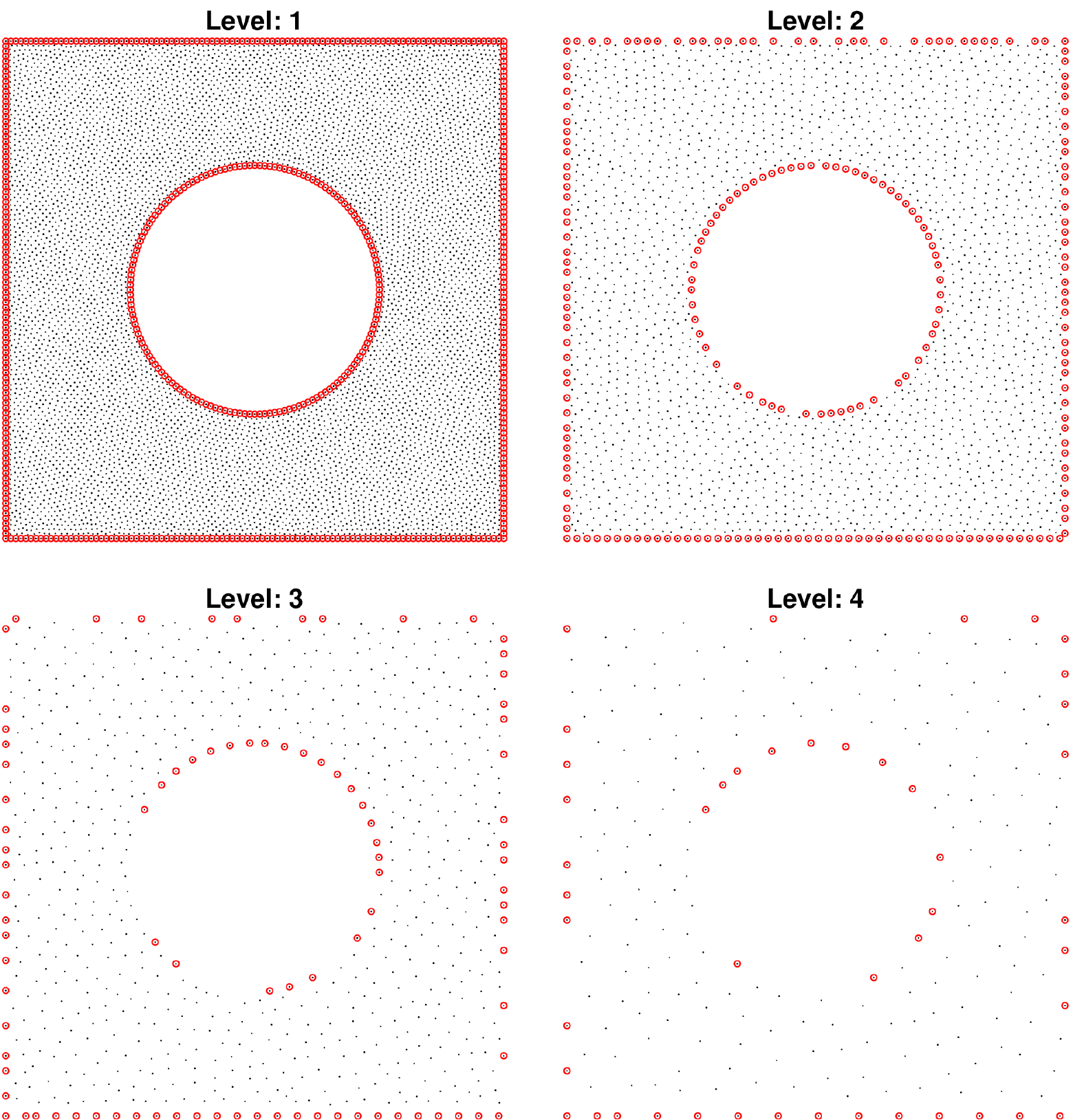}
        \caption{Naive algorithm \label{fig:bdry1}}
    \end{subfigure}
    \begin{subfigure}[t]{0.45\textwidth}
        \includegraphics[width=\textwidth]{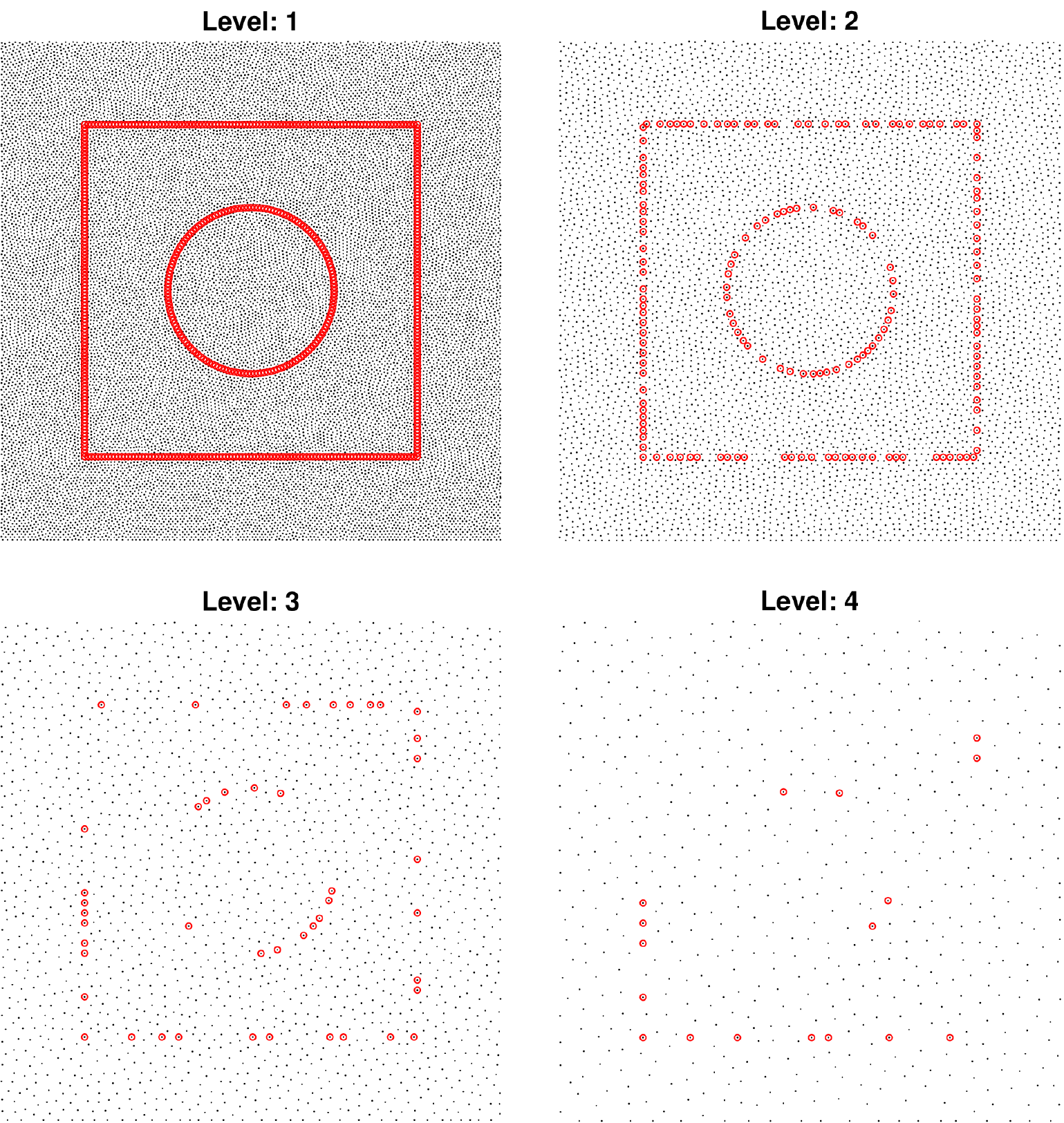}
        \caption{Nodes included interior and exterior to boundaries \label{fig:bdry2}}
    \end{subfigure}
    \caption{The moving front subsampling algorithm applied to a test node set with two boundaries. The initial node set is subsampled three times. The boundary node set is included in the domain node set such that they are subsampled collectively and simultaneously. Subsampling performance along the boundary is improved by including nodes interior and exterior to all boundaries.}
\end{figure}

When applying the moving front algorithm naively to a set for which the boundary nodes are included in the domain, the top boundary nodes are undesirably subsampled faster than the ones at the lower boundary, as demonstrated in Figure \ref{fig:bdry1}. One way to reduce the subsampling inconsistencies is to include nodes interior and exterior to each boundary as seen in Figure \ref{fig:bdry2}.

Another way to reduce the inconsistencies in subsampling which may be due to the inherent directional bias of the moving front algorithm is to alternate the direction between subsampling iterations. This technique of alternating direction is unsatisfactory, however, because an ideal method would be effective independent any inherent directional bias. To further improve robustness of the moving front method in the presence of boundaries, the following section considers subsampling boundaries separately.

%%%%%%%%%%%%%%%%%%%%%%%%%%%%%
\subsection{Subsample Boundary Separately} \label{sec:subBoundSep}
\begin{figure}[]
    \centering
    \begin{subfigure}[t]{0.45\textwidth}
        \includegraphics[width=\textwidth]{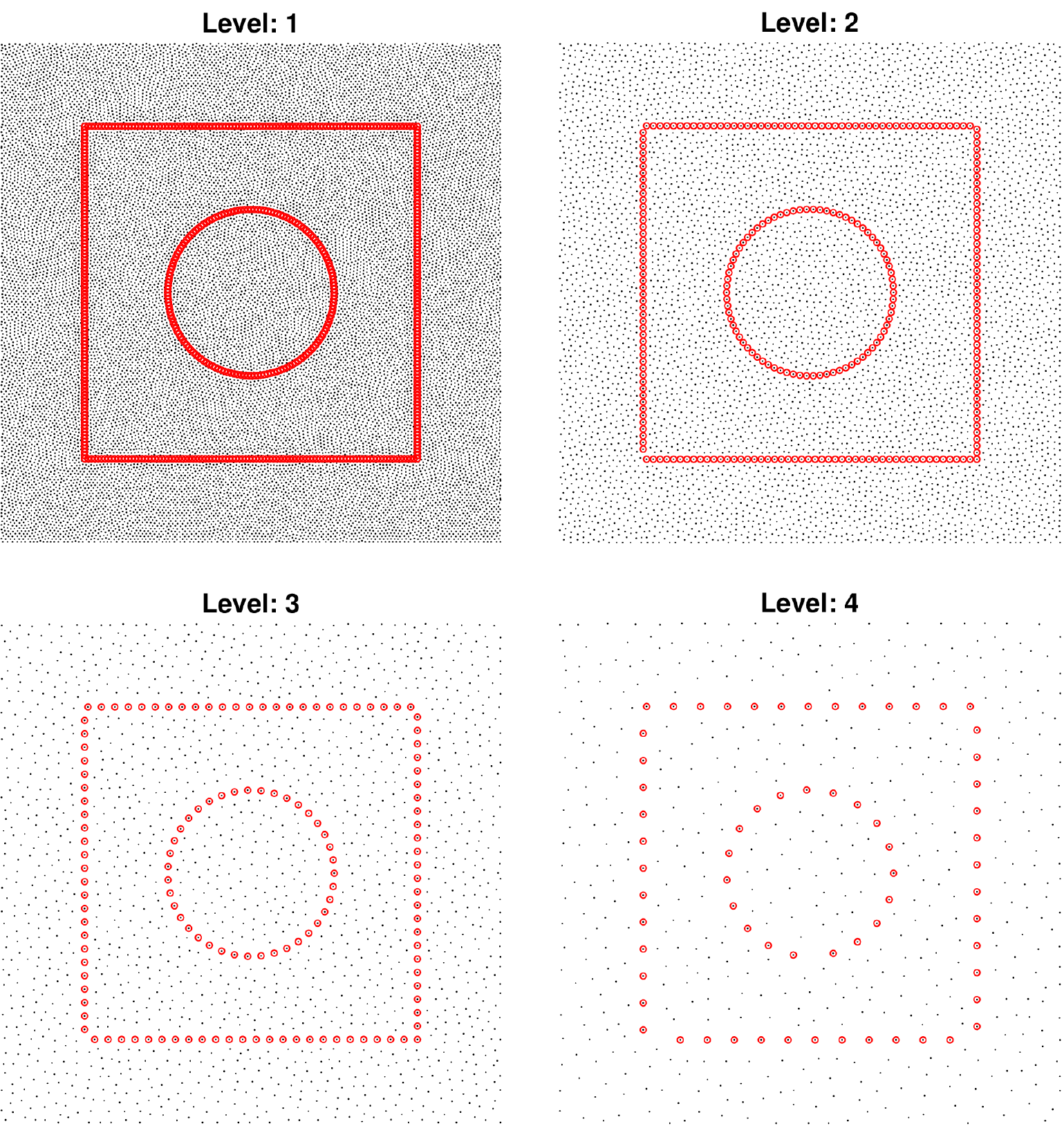}
        \caption{Alternating direction \label{fig:bdry4}}
    \end{subfigure}
    \begin{subfigure}[t]{0.45\textwidth}
        \includegraphics[width=\textwidth]{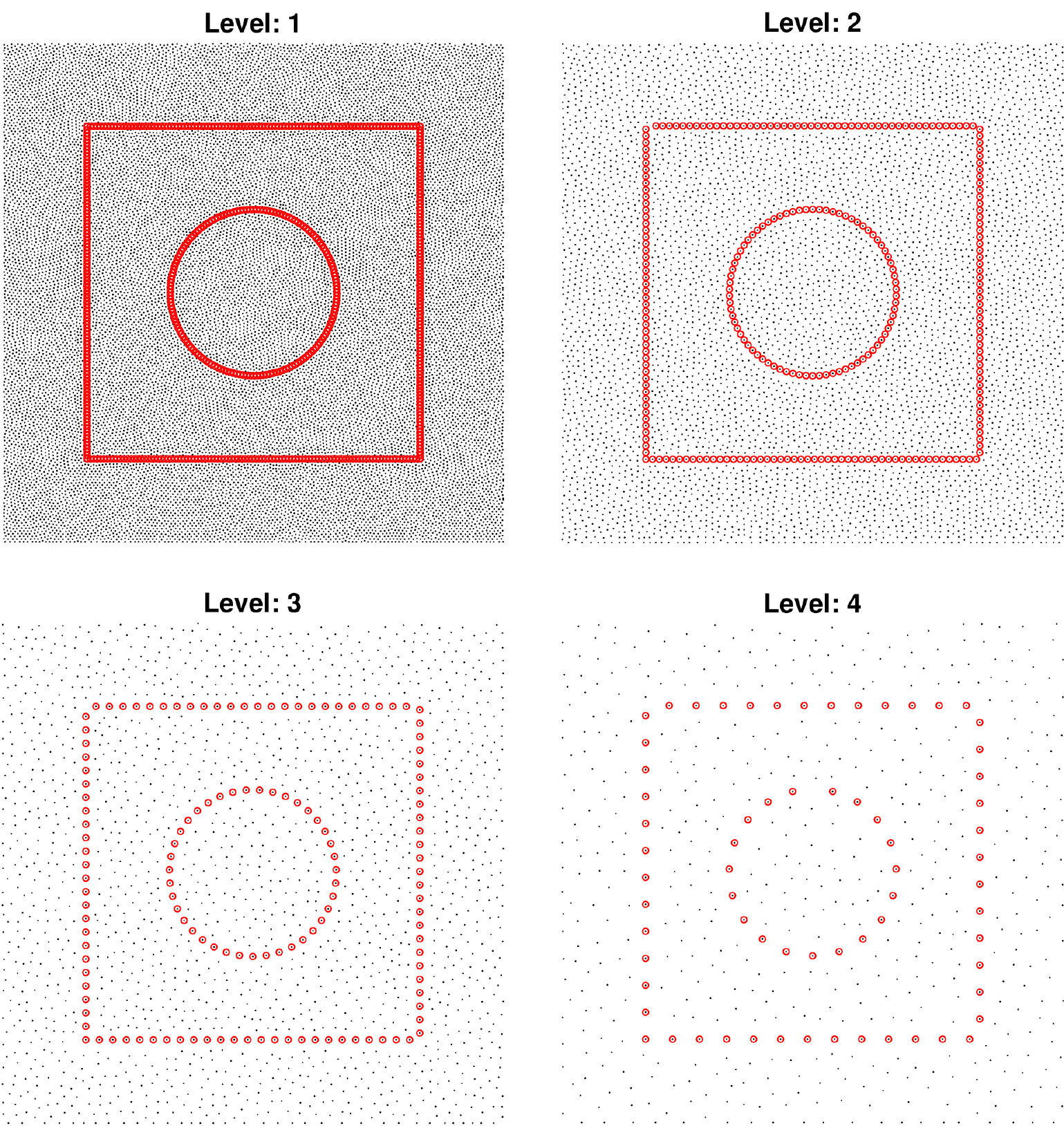}
        \caption{No alternating direction \label{fig:bdry5}}
    \end{subfigure}
    \caption{The moving front subsampling algorithm applied to a test node set with two boundaries. The initial node set is subsampled three times. In these figures, the boundary node set is subsampled separately from the domain node set. Subsampling performance along the boundary is improved by subsampling boundary nodes independently. Additionally, no directional bias is detectable even when the algorithm does not alternate direction. \label{fig:bdry45}}
\end{figure}

First, the given boundary nodes are subsampled. Then, any domain nodes within a prescribed distance of the boundary nodes are removed. Finally, the domain nodes are subsampled. The effects of this two-step process can be seen in Figure \ref{fig:bdry45}. Figure \ref{fig:bdry4} alternates direction of the moving front algorithm while Figure \ref{fig:bdry5} does not.  A significant improvement in how consistently the algorithm behaves across the node set is apparent, independent of any direction bias in the subsampling algorithm.

%%%%%%%%%%%%%%%%%%%%%%%%%%%%%%%%%%%%%%%%%%%%%%%%%%%%%%%%%%
\section{Comparisons of Subsampling Methods} \label{sec:Tests}
\begin{figure}[H]
    \centering
    \begin{subfigure}[c]{0.4\textwidth}
        \includegraphics[bb=0 0 256 256, width=\textwidth]{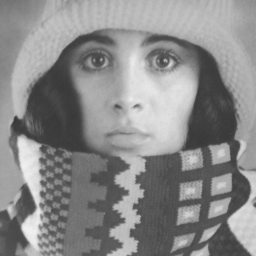}
        \caption{Original trui.png image \label{fig:truiOriginal}}
    \end{subfigure}
    \begin{subfigure}[c]{0.4\textwidth}
        \includegraphics[width=\textwidth]{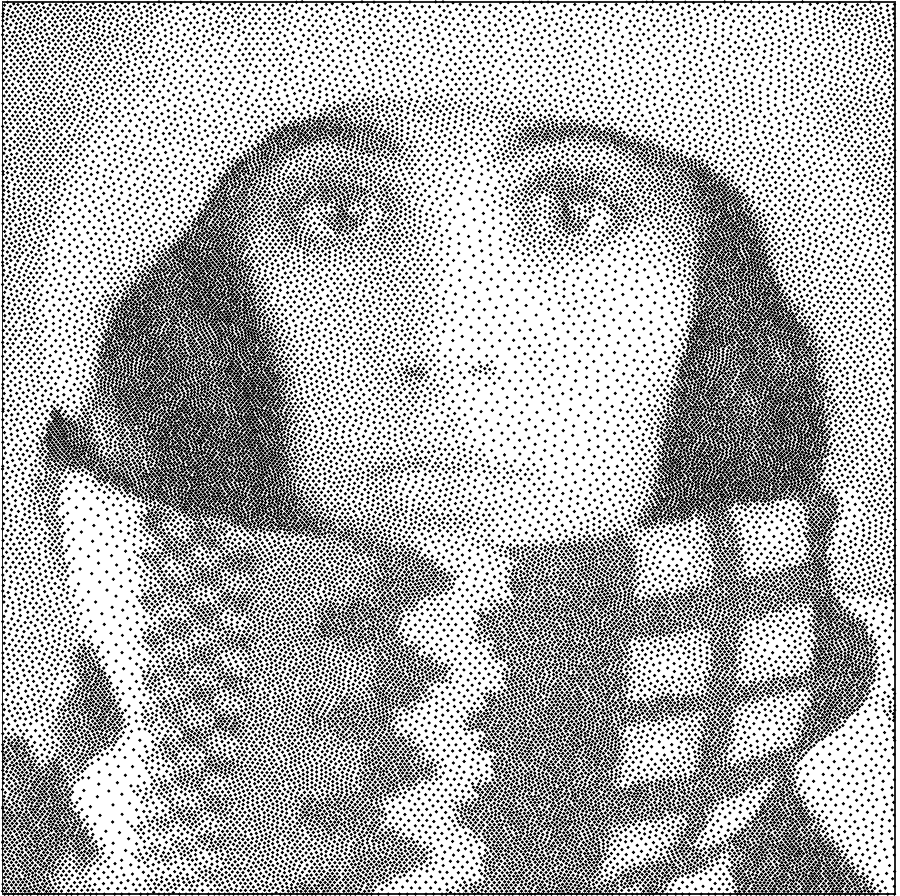}
        \caption{Dithered trui image \label{fig:truiDithered}}
    \end{subfigure}
    \caption{The original trui.png image and a dithered version with 36,303 nodes, obtained by the algorithm in \cite{vanderSande2021}.}
\end{figure}

This section compares the performance in preserving node density variation through iterative coarsening of the four subsampling algorithms found in Section \ref{sec:SubAlgs}. For each example, the primary node set has been generated from the trui image, see \Cref{fig:truiOriginal}, using the node generation algorithm from \cite{vanderSande2021}, see \Cref{fig:truiDithered}. While this initial node set does not have immediate application to solving PDEs, its radically varying node densities make the trui image a good test problem for visually spotting any algorithmic artifacts. The trui image also contains regions of uniform density, thus illustrating subsampling capabilities on regions of both locally variable and locally uniform node densities.

%%%%%%%%%%%%%%%%%%%%%%%%%%%%%
\subsection{Heuristic Comparison}
\begin{figure}[]
    \centering
    % --- Moving Front
    \rotatebox[origin=c]{90}{Moving Front}\quad
    \begin{subfigure}[c]{0.3\textwidth}
        \includegraphics[width=\textwidth]{Figures/trui_dithered.pdf}
        \label{fig:truiMF}
    \end{subfigure}
    \begin{subfigure}[c]{0.3\textwidth}
        \includegraphics[width=\textwidth]{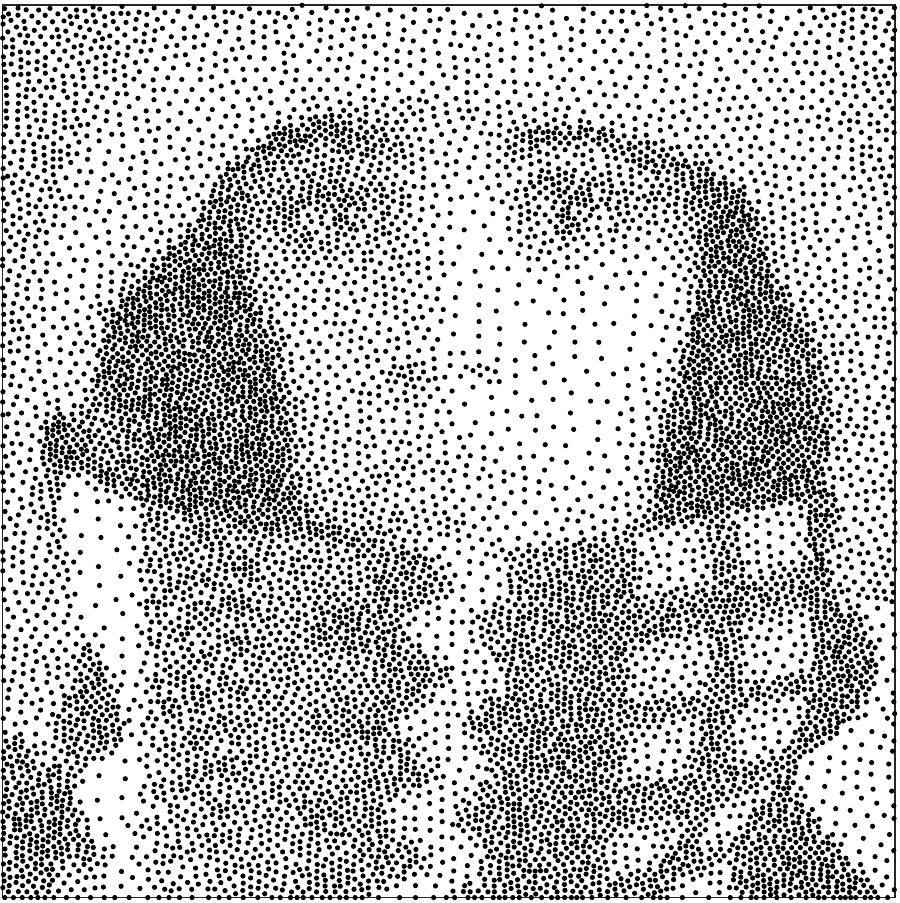}
        \label{fig:truiSub1MF}
    \end{subfigure}
    \begin{subfigure}[c]{0.3\textwidth}
        \includegraphics[width=\textwidth]{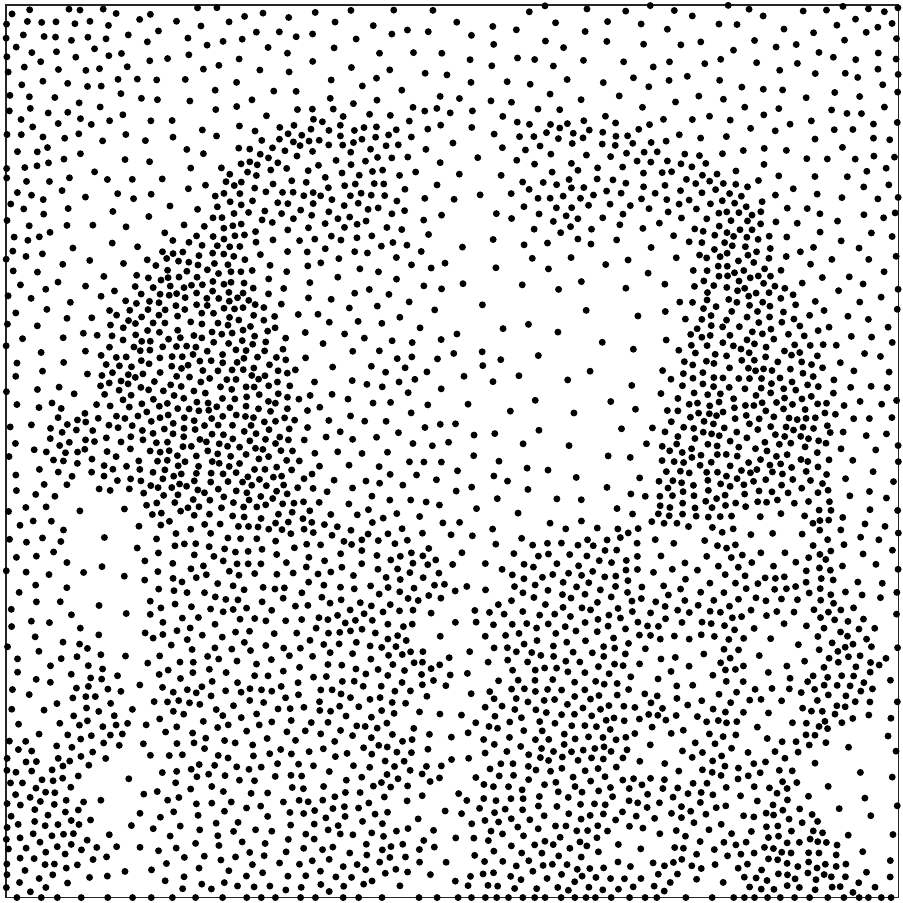}
        \label{fig:truiSub2MF}
    \end{subfigure}
    
    % --- Weighted
    \rotatebox[origin=c]{90}{Weighted}\quad
    \begin{subfigure}[c]{0.3\textwidth}
        \includegraphics[width=\textwidth]{Figures/trui_dithered.pdf}
        \label{fig:truiWS}
    \end{subfigure}
    \begin{subfigure}[c]{0.3\textwidth}
        \includegraphics[width=\textwidth]{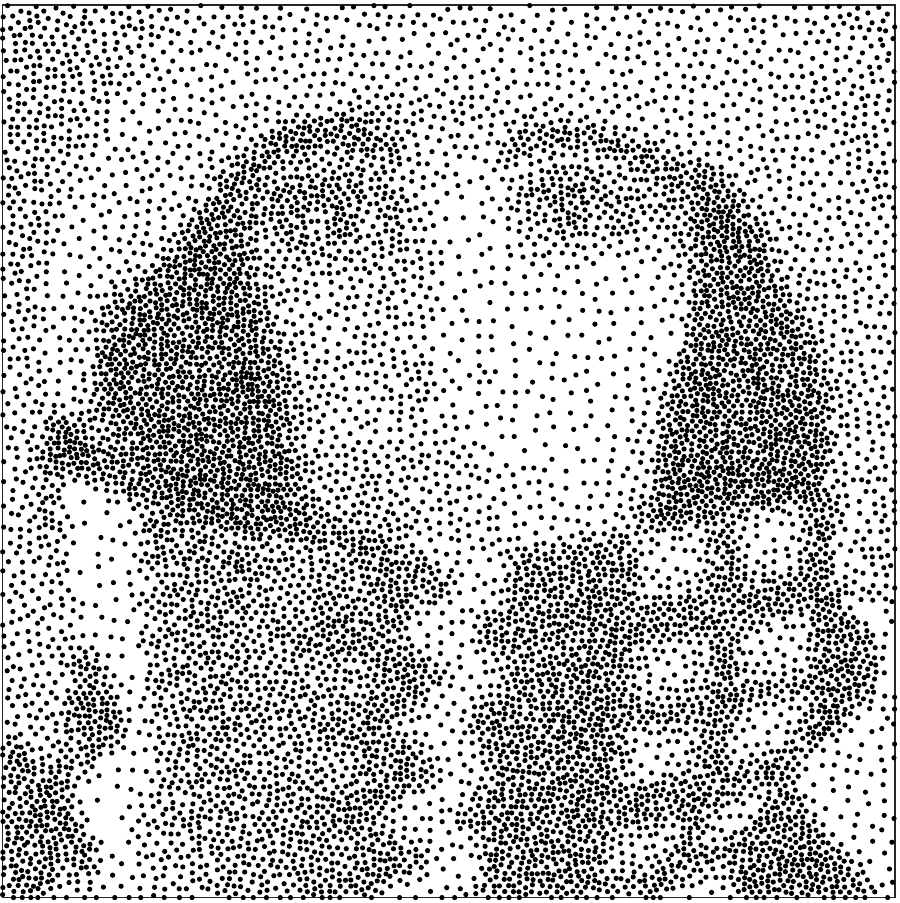}
        \label{fig:truiSub1WS}
    \end{subfigure}
    \begin{subfigure}[c]{0.3\textwidth}
        \includegraphics[width=\textwidth]{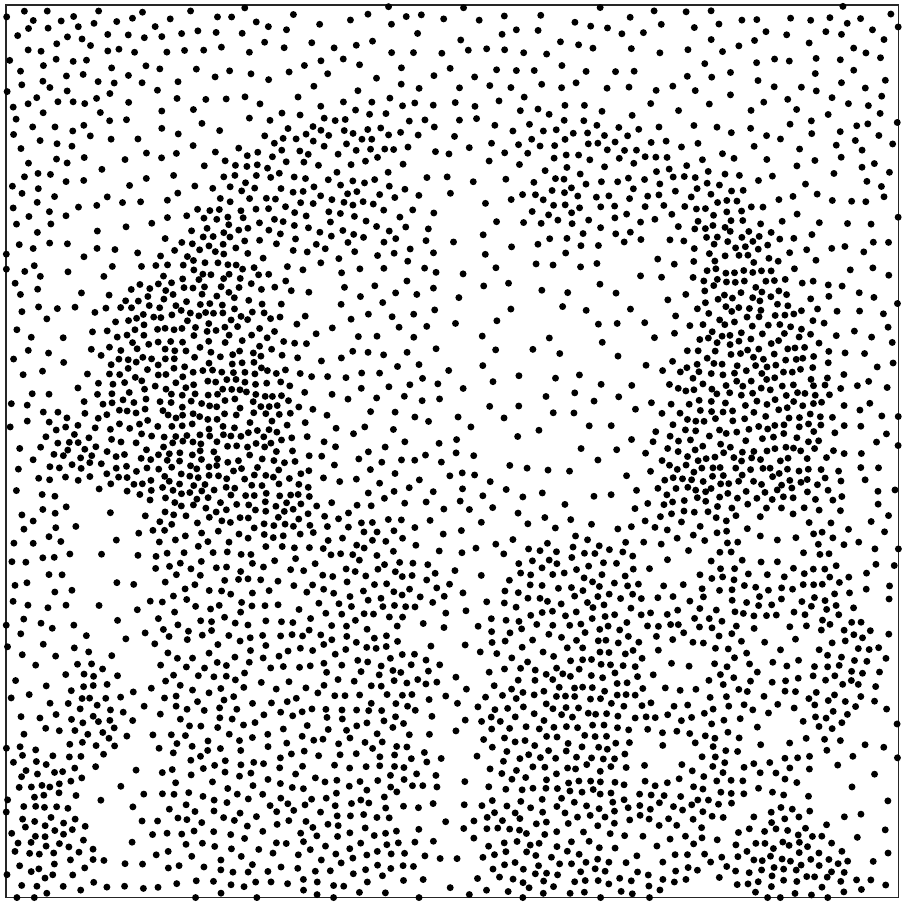}
        \label{fig:truiSub2WS}
    \end{subfigure}
    
    % --- Poisson Disk
    \rotatebox[origin=c]{90}{Poisson Disk}\quad
    \begin{subfigure}[c]{0.3\textwidth}
        \includegraphics[width=\textwidth]{Figures/trui_dithered.pdf}
        \label{fig:truiPD}
    \end{subfigure}
    \begin{subfigure}[c]{0.3\textwidth}
        \includegraphics[width=\textwidth]{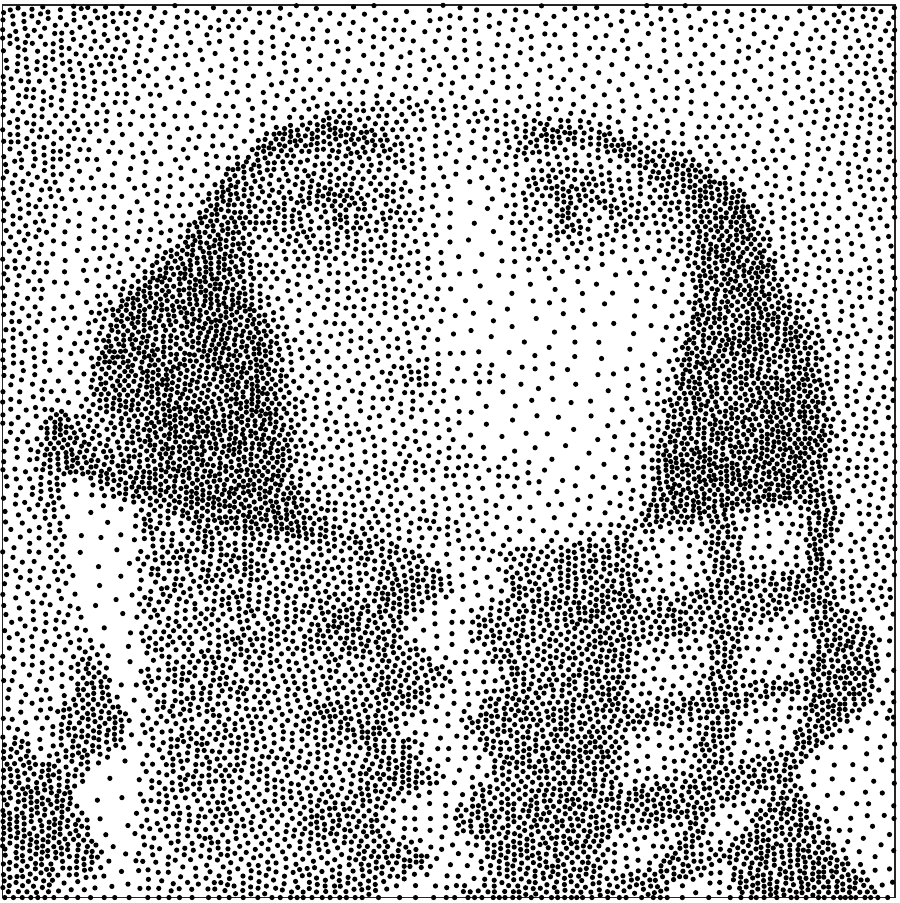}
        \label{fig:truiSub1PD}
    \end{subfigure}
    \begin{subfigure}[c]{0.3\textwidth}
        \includegraphics[width=\textwidth]{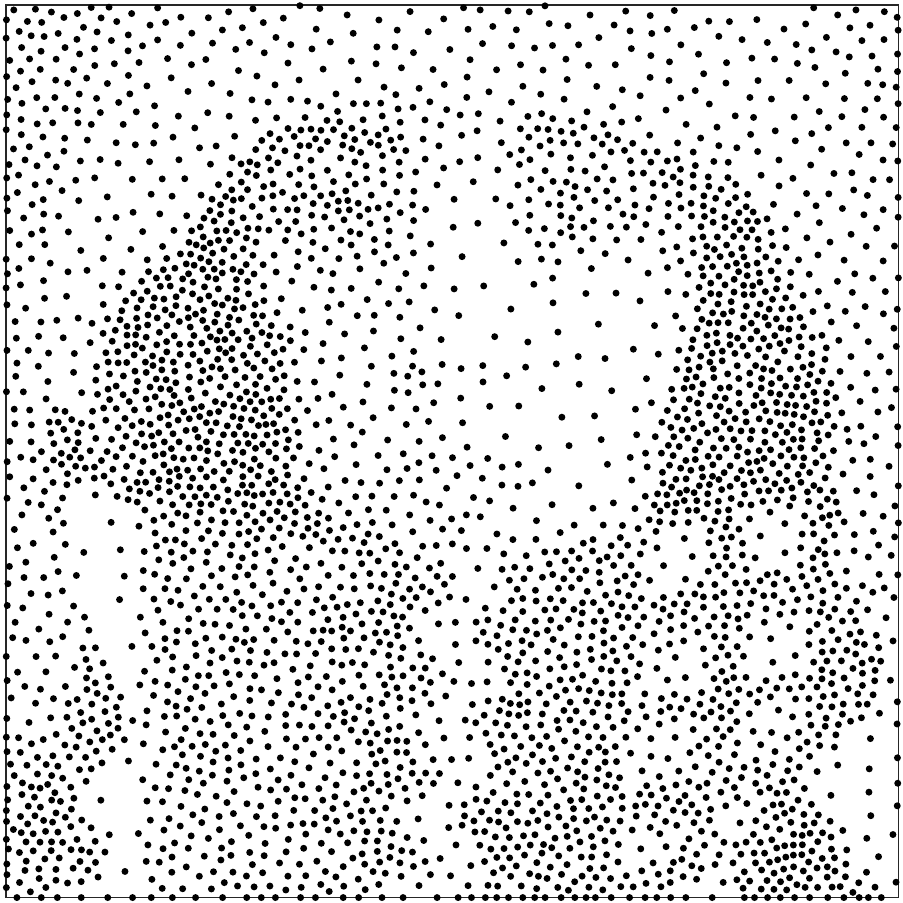}
        \label{fig:truiSub2PD}
    \end{subfigure}
    
    % --- generalized Diversity
    \rotatebox[origin=c]{90}{Generalized Diversity}\quad
    \begin{subfigure}[c]{0.3\textwidth}
        \includegraphics[width=\textwidth]{Figures/trui_dithered.pdf}
        \label{fig:truigDS}
    \end{subfigure}
    \begin{subfigure}[c]{0.3\textwidth}
        \includegraphics[width=\textwidth]{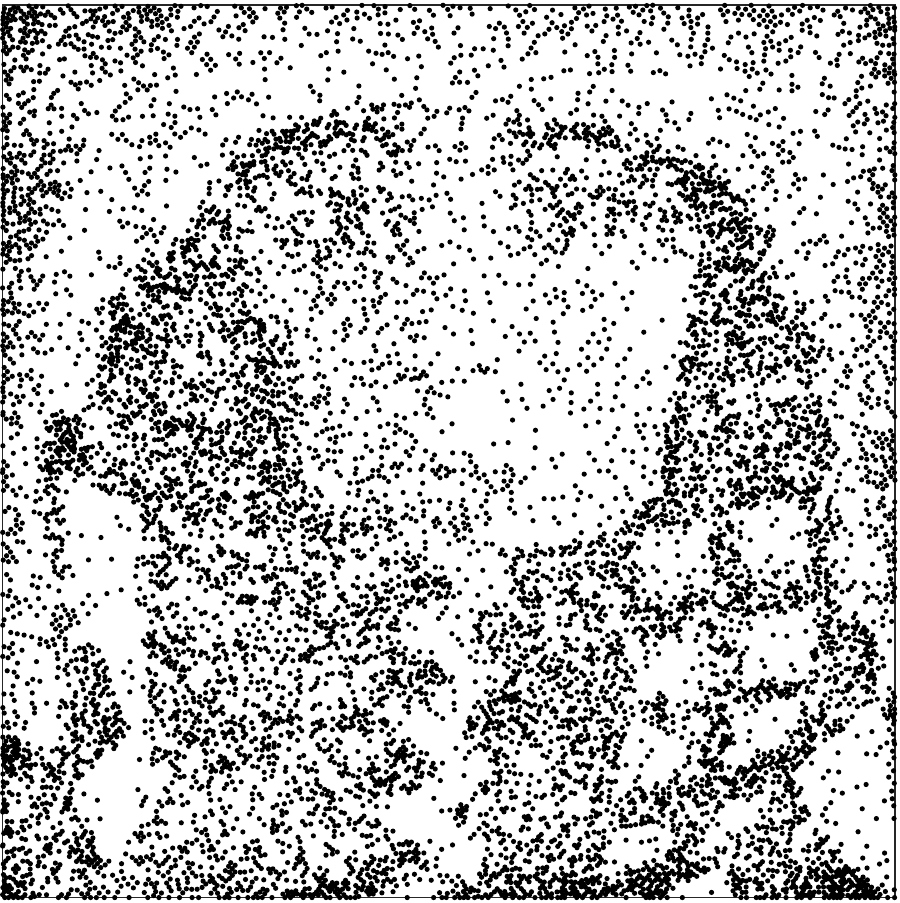}
        \label{fig:truiSub1gDS}
    \end{subfigure}
    \begin{subfigure}[c]{0.3\textwidth}
        \includegraphics[width=\textwidth]{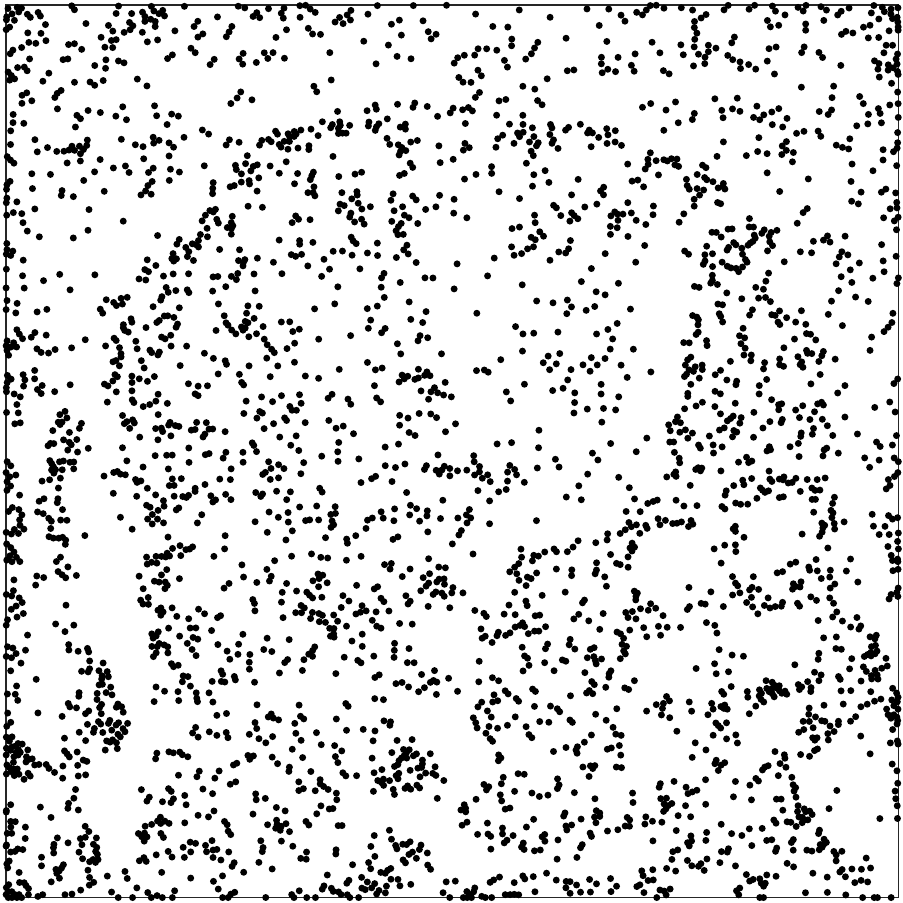}
        \label{fig:truiSub2gDS}
    \end{subfigure}
    \caption{A visualization of the initial dithered trui image (36303 nodes), the first subsampling (10553 nodes), and the second subsampling (3404 nodes). From top to bottom, the rows demonstrate the moving front, weighted, Poisson disk, and generalized diversity subsampling algorithms. The first column of the figure contains a redundant image of the initial dithering for convenience in comparison. \label{fig:truiMatrix}}
\end{figure}

This section provides visualizations of the subsampled node sets of each algorithm. Each algorithm is applied twice to the dithered trui image as seen in Figure \ref{fig:truiMatrix}. The original node set contains 36303 nodes, the first subsample contains 10553 nodes, and the second subsample contains 3404 nodes. Each algorithm discussed, excluding the generalized diversity subsampling algorithm\footnote{The generalized diversity subsampling algorithm explicitly requires a target number of nodes as input rather than a parameter.}, requires a parameter, $c$, to control the level of coarsening. To reproduce the subsamples in Figure \ref{fig:truiMatrix}, the values of $c$ used in each algorithm are listed in Table \ref{tbl:params}. The moving front algorithm also relies on a choice of nearest neighbors which was $k=10$ for these tests.

\begin{table}[H]
    \centering
    \begin{tabular}{c|c|c}
        Method & First Subsampling & Second Subsampling\\
        \hline
        MF & 1.5101 & 1.518 \\
        W & 3.44 & 3.1 \\
        PD & 1.4931 & 1.5394
    \end{tabular}
    \caption{The parameters, $c$, for reproducing the node sets in Figure \ref{fig:truiMatrix} for the moving front (MF), weighted (W), and Poisson disk (PD) subsampling algorithms. The generalized diversity subsampling algorithm explicitly relies on the desired number of nodes in the coarse node set and thus has no parameter listed here. The moving front algorithm also relies on a choice of nearest neighbors which was $k=10$ for these tests.}
    \label{tbl:params}
\end{table}

A heuristic comparison between the subsamplings iterations primarily demonstrates the visual goodness of the first three algorithms over the generalized diversity subsampling algorithm. Among the remaining three, the woman's nostrils are more distinct in the moving front and Poisson disk algorithms than in the weighted subsampling while the Poisson disk algorithm seems to preserve the mouth slightly more clearly by the second subsampling. Additionally, the moving front and Poisson disk algorithms preserve a higher level of clarity in the patterns\footnote{The checkered pattern in the scarf on the woman's right is a clear example.} in the trui scarf than the weighted subsampling algorithm. Again, the moving front and Poisson disk subsampling algorithms each better preserve the density disparity between areas of low and high node density in the original dithering\footnote{The disparity is most notable between the woman's left cheek and hair and between the light stripe in the scarf on the woman's right and any of the surrounding regions.} than do the weighted or generalized diversity subsampling algorithms. Finally, the Poisson disk algorithm may have a tendency to subsample too aggressively in places\footnote{The light stripe on the scarf's right side (left side of the image) is much sparser than in the moving front or weighted subsampling algorithms}. It should be noted that no direction bias of the moving front algorithm is visible.

%%%%%%%%%%%%%%%%%%%%%%%%%%%%%
\subsection{Node Quality Measures}\label{sec:nodeQuality}

While visually comparing the results of each algorithm is useful, ultimately a more rigorous and quantitative comparison is desirable. As such, a variety of node quality measures are commonly discussed when comparing uniform subsampling methods \cite{vanderSande2021,Shang2022DiversitySC}. However, metrics which describe the quality of spatially variable node densities are less common \cite{vanderSande2021}. One way of evaluating the quality of a variable density node set on its own is through local regularity distribution of distance to the nearest $k$ neighbors $\delta_{i,j}$, $i=1,2,...k$ for each node $x_j$.

It is here considered sufficient to expect the node generation algorithm to produce an initial node set which is sufficiently good\footnote{Where goodness is  determined by any number of node quality measures chosen based upon context.}. The responsibility of a subsampling algorithm is then to preserve the characteristics of the original node set. A natural way to determine how well any subsampling preserves the quality of the initial node set is to measure the coarse node set in comparison to the fine node set. 

The two novel measures presented here are straightforward extensions of the commonly accepted measures of local regularity for evaluating the quality of variable density node sets and are referred to as measures of comparative local regularity (CLR). These CLRs contrast the average or standard deviation of distances to the $k$ nearest neighbors of the fine node set from those of the coarse node set. For each node set $X$, the Euclidian distance between each node $x_j\in X$ and its $k$ nearest neighbors is calculated as $\delta_{i,j} = \|x_j - x_{i,j}\|$. The average $\overline{\delta}_j$ or standard deviation $\sigma_j$ of these distances is then found for each node $x_j in X$. These are typical measures of local regularity. However, the present goal is measure how similar the initial and subsampled node sets are. To this end, the distributions of $\overline{\delta}_j^{\ \text{fine}}$ and $\overline{\delta}_j^{\ \text{coarse}}$ over each of the fine and coarse node sets, respectively, are first normalized to be between 0 and 1. Then the difference between these distributions is calculated at each of the nodes in the coarse node set and each of the collocated nodes in the fine node set\footnote{Effectively, $\overline{\delta}_j^{\text{fine}}$ needs only be calculated at those nodes $x_j$ in the fine node set which also belong to the coarse node set}. Finally, the standard $L_2$ norm is taken to produce a measure of CLR. The same process can be applied using the standard deviation $\sigma_j$ of those distances $\delta_{i,j}$.

\begin{figure}[]
    \centering
    \begin{subfigure}[c]{0.41\textwidth}
        \includegraphics[width=\textwidth]{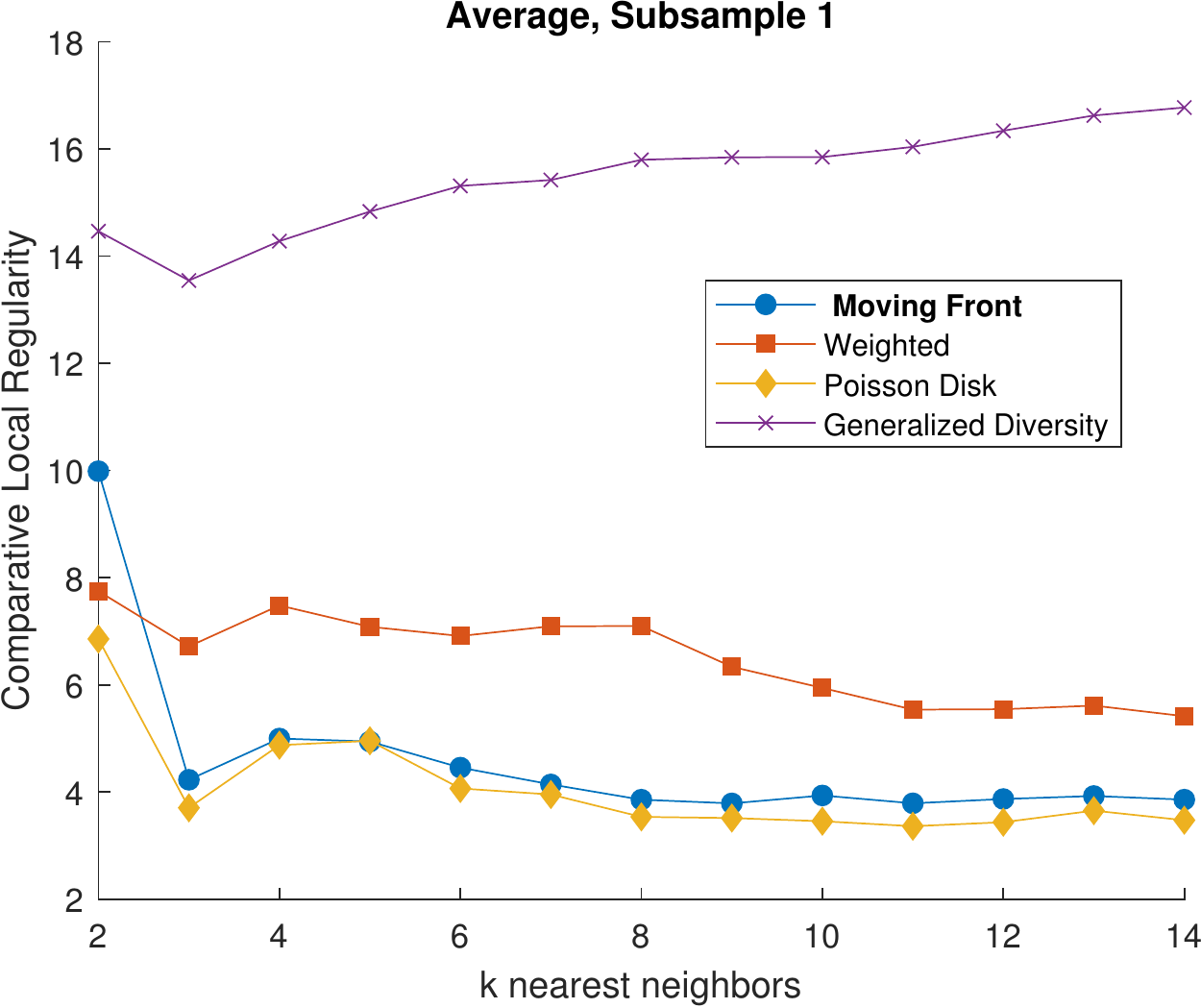}
        \label{fig:CLRASub1}
    \end{subfigure}
    \begin{subfigure}[c]{0.41\textwidth}
        \includegraphics[width=\textwidth]{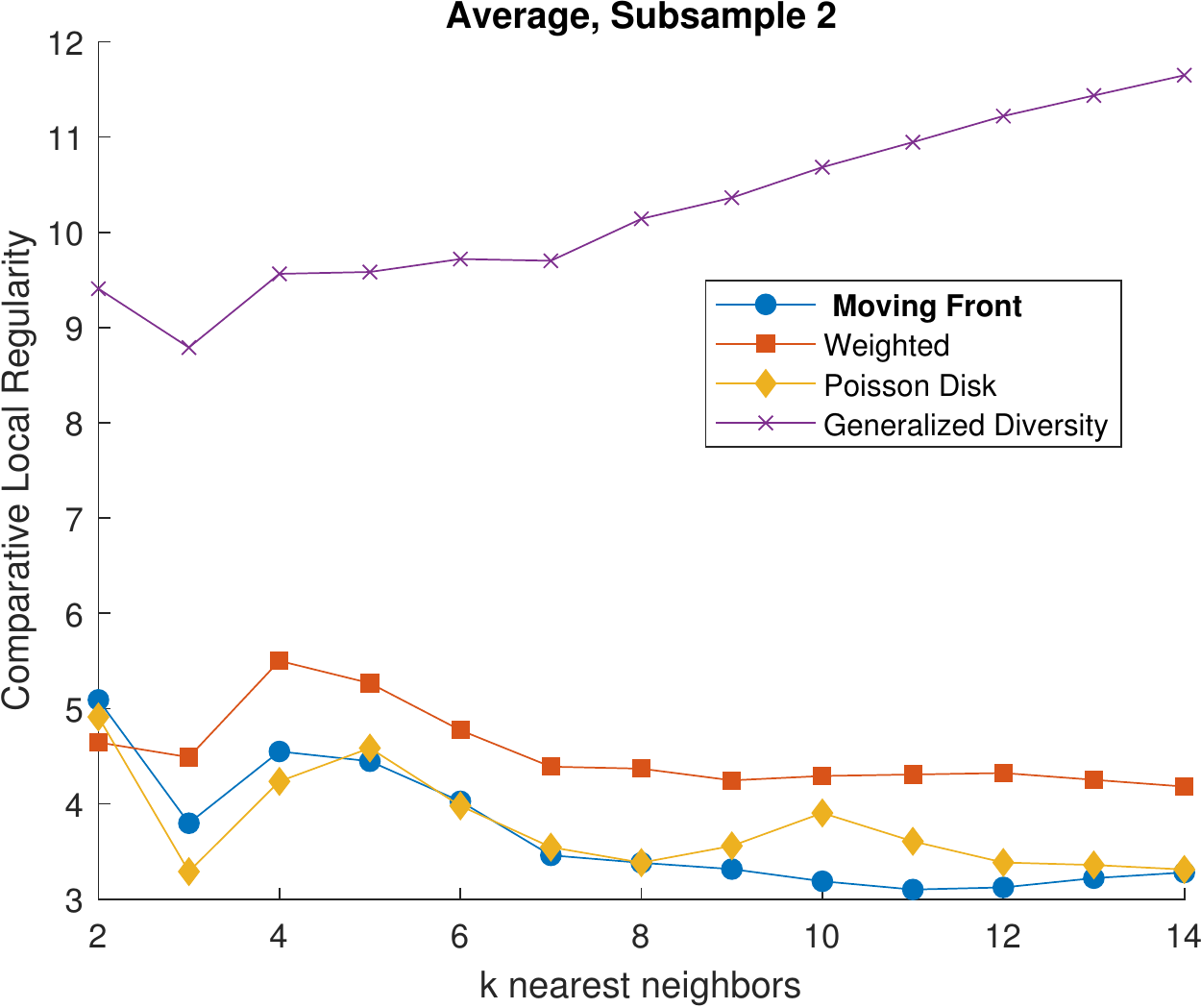}
        \label{fig:CLRASub2}
    \end{subfigure}
    \begin{subfigure}[c]{0.41\textwidth}
        \includegraphics[width=\textwidth]{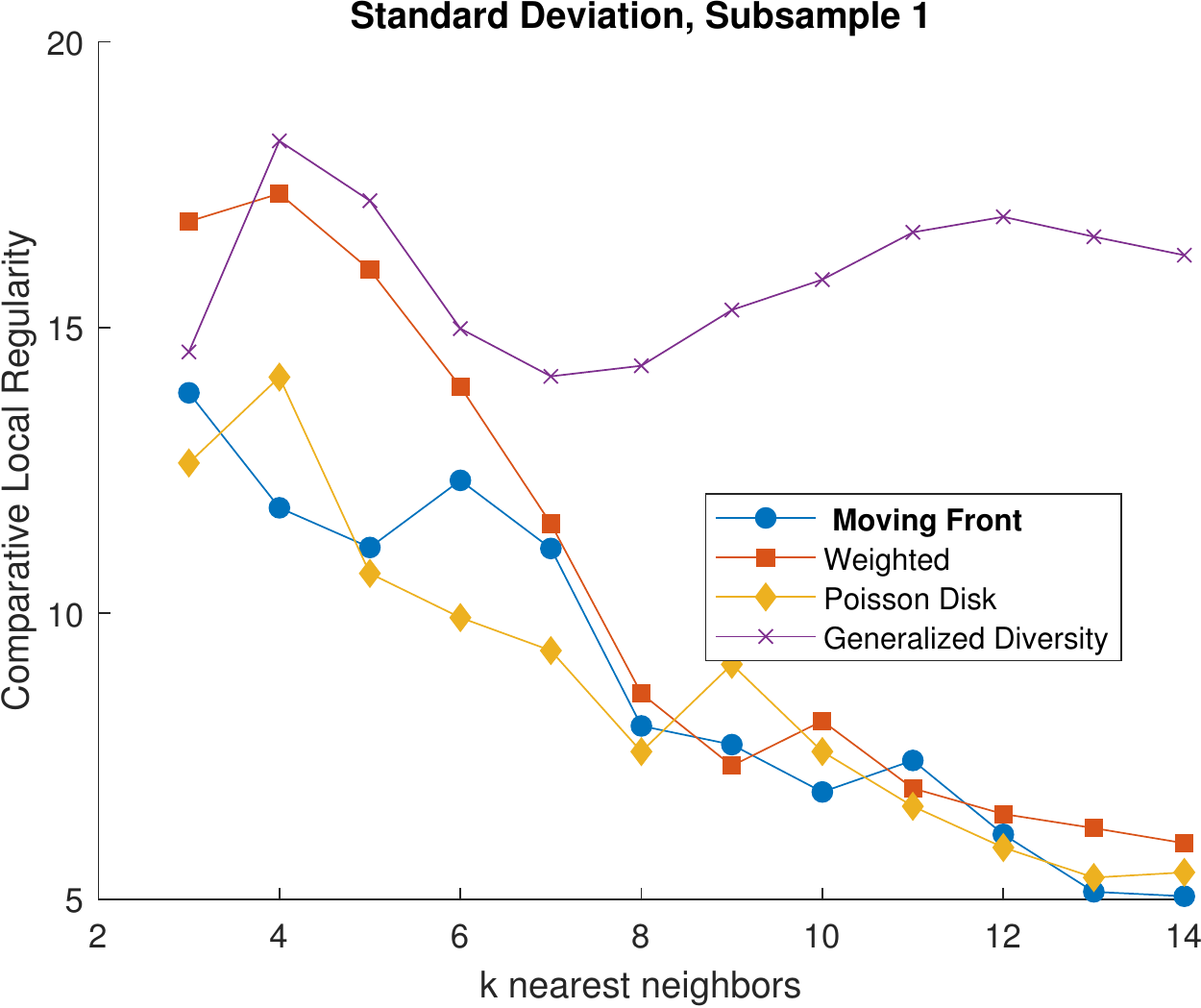}
        \label{fig:CLRSDSub1}
    \end{subfigure}
    \begin{subfigure}[c]{0.41\textwidth}
        \includegraphics[width=\textwidth]{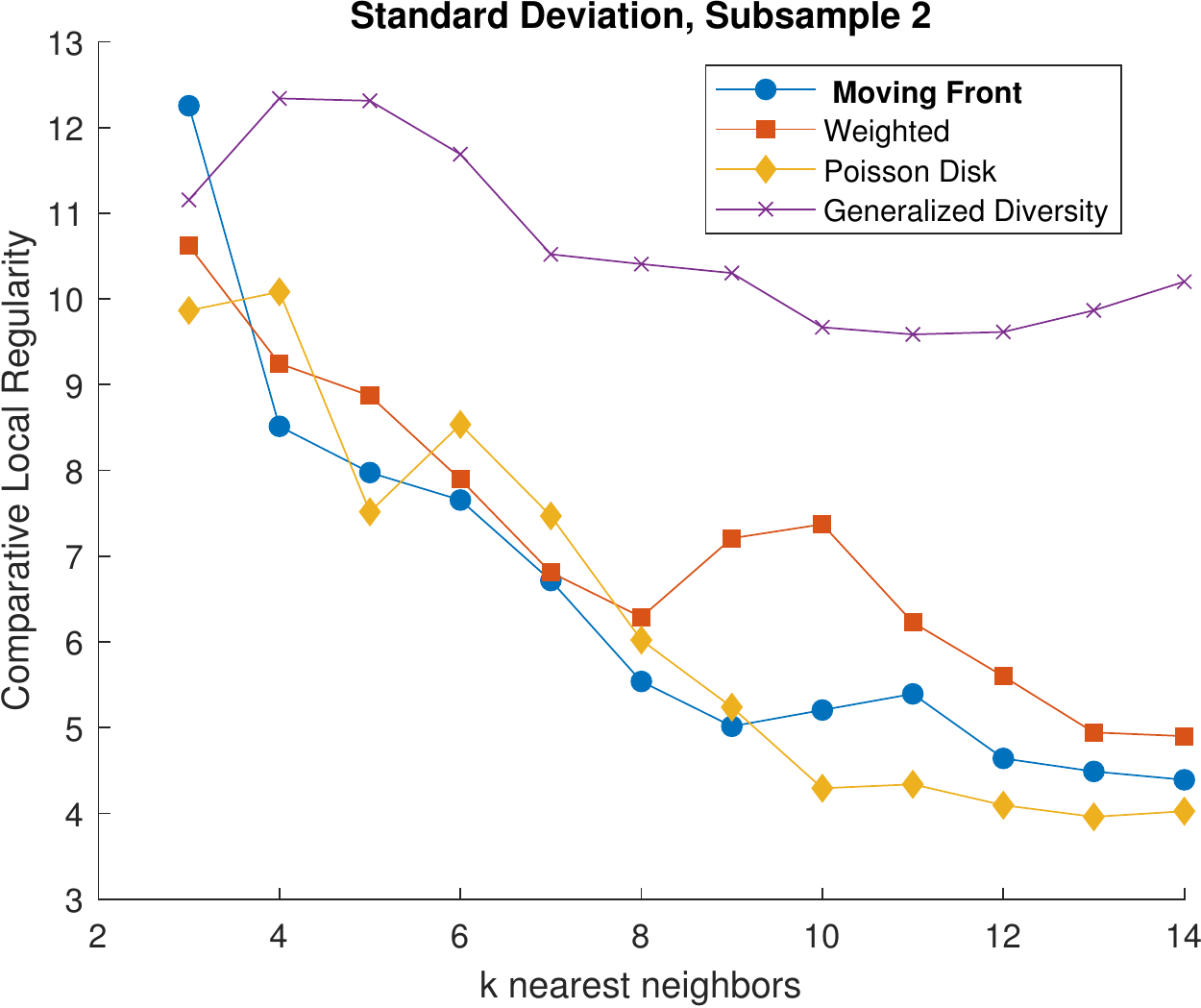}
        \label{fig:CLRSDSub2}
    \end{subfigure}
    \caption{The comparative local regularity (CLR) of the average distance and standard deviation of distances for $k = 2,3,...,14$ nearest neighbors of various subsampling methods (weighted, moving front, Poisson disk, and generalized diversity subsampling) applied once and twice to the dithered trui image. For both measures of CLR, a lower value is better.}
    \label{fig:CLR}
\end{figure}

Given that these CLRs measure the difference between a given initial node set and a subsampling of it, it is ideal to minimize these measures across algorithms. The presented CLRs should only be compared between subsamplings of the same size and from the same initial node set. Upon comparison of the average and standard deviation CLRs for each algorithm across various values of nearest neighbors as seen in Figure \ref{fig:CLR}, it is clear that the moving front and Poisson disk algorithms preserve the characteristics of the initial dithered trui image better than the other two algorithms. This behavior is consistent across various sizes of node sets. Between them, however, it is not clear which one is qualitatively better.

%%%%%%%%%%%%%%%%%%%%%%%%%%%%%
\subsection{Computational Cost} \label{sec:CompCost}
\begin{figure}
	\centering
	\includegraphics[width=0.65\textwidth]{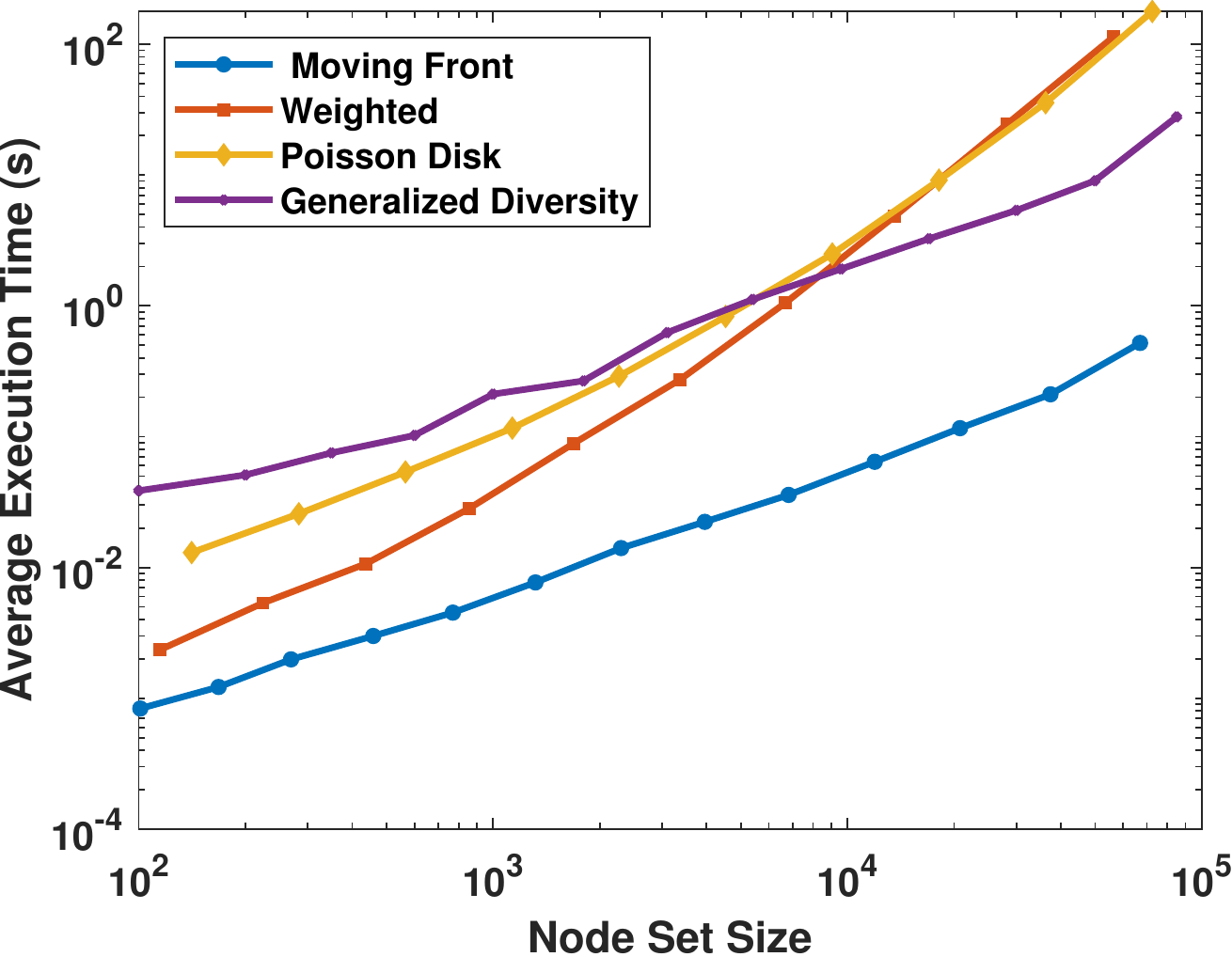}
	\caption{The computation time of each subsampling algorithm. The node set size is that of resultant subsample. Each subsample is subsequent such that the coarse node set of the previous iteration is the fine node set of the next. The execution time was averaged over ten iterations of the moving front, weighted, Poisson disk, and generalized diversity subsampling algorithms. Note the logarithmic scales. The moving front algorithm is significantly faster than any of the other algorithms.}
	\label{fig:CompCost}
\end{figure}

When evaluating any numerical scheme, computational cost must be considered as much as any other aspect of performance. The moving front, weighted, and Poisson disk subsampling algorithms presented in this paper were coded in MATLAB\footnote{The code for the moving front algorithm included in Appendix \ref{apx:MF} is provided in Python for the readers convenience. Code is available in both MATLAB and Python on the authors GitHub \cite{Lawrence_GitHub}.} The generalized diversity subsampling algorithm is coded in Python as presented in \cite{Shang2022DiversitySC}. A comparison of the computational complexity can be found in Figure \ref{fig:CompCost}. The average execution times in Figure \ref{fig:CompCost} are calculated based on ten repetitions of each algorithm subsampling from the node set size of the previous data point to the node set size for the given point using the \verb|timeit| commands native to MATLAB and Python. The moving front algorithm is clearly the fastest algorithm of those compared here. The significant computational cost savings secures the moving front algorithm as the best overall performing algorithm.

%\begin{table}[H]
%    \centering
%    \begin{tabular}{c|c|c}
%        Method & Mean Run Time (s) & Run Time Standard Deviation (ms)\\
%        \hline
%        MF & 0.0866 & 2.17 \\
%        W & 22.8 & 965 \\
%        PD & 3.06 & 62.9 \\
%        gDS & 2.55 & 206 \\
%    \end{tabular}
%    \caption{The computation time of the first subsampling of the 36303 node dithered trui image to 10553 nodes for various methods (MF: moving front subsampling, W: weighted subsampling, PD: Poisson disk subsampling, and gDS: generalized diversity subsampling). For each algorithm, the first subsampling is repeated ten times to calculate the mean and standard deviation. The moving front algorithm is two to three magnitudes faster than any other algorithm.}
%    \label{tbl:CompCost}
%\end{table}

%%%%%%%%%%%%%%%%%%%%%%%%%%%%%%%%%%%%%%%%%%%%%%%%%%%%%%%%%%
\section{Meshfree multilevel RBF-FD solver} \label{sec:Example}
\begin{figure}[H]
    \centering
    \begin{subfigure}[c]{0.4\textwidth}
        \includegraphics[width=\textwidth]{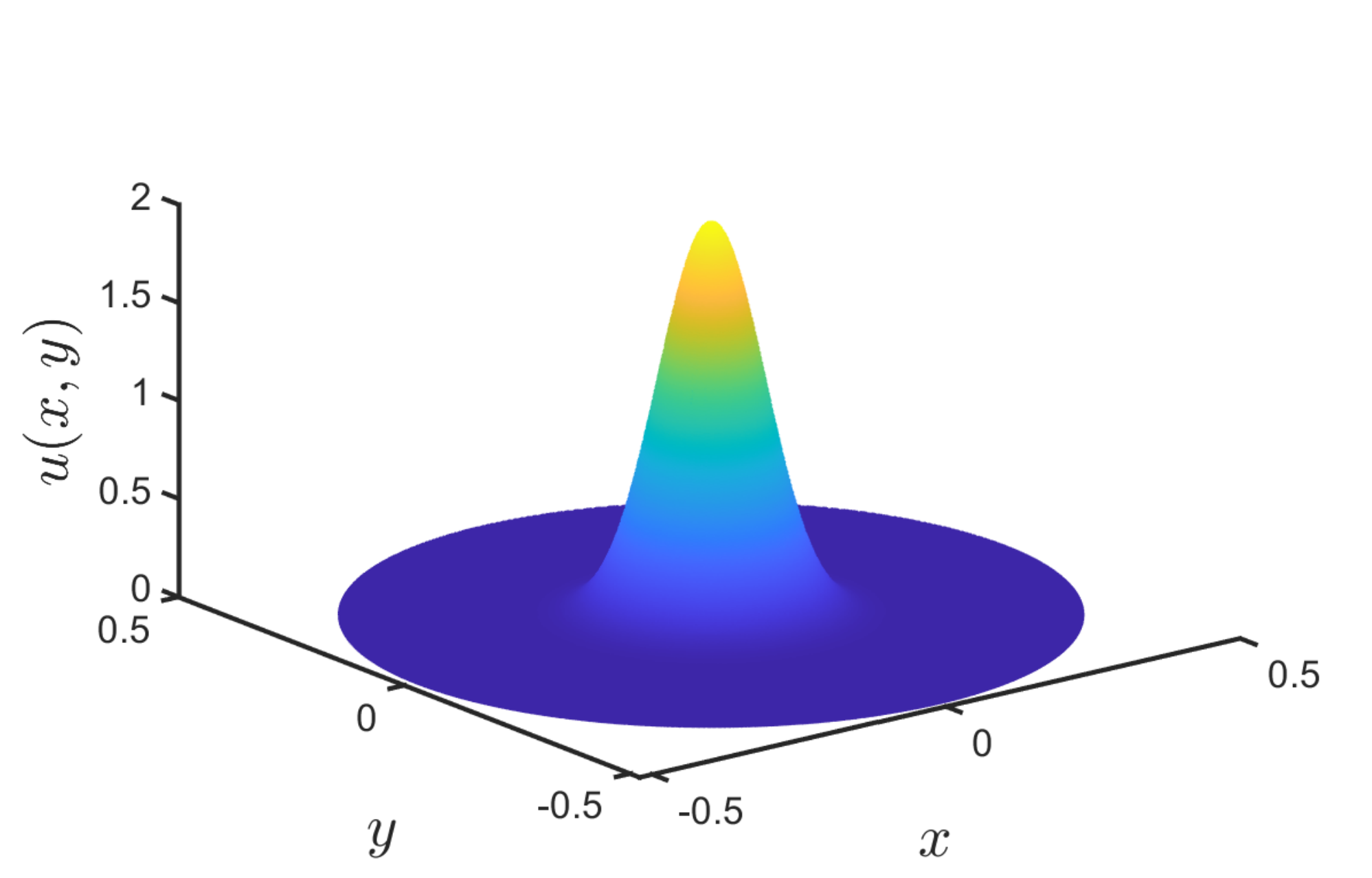}
        \caption{$u(r,\theta) = 2\exp(-r/ 0.01)$ }
    \end{subfigure}
    \begin{subfigure}[c]{0.4\textwidth}
        \includegraphics[width=\textwidth]{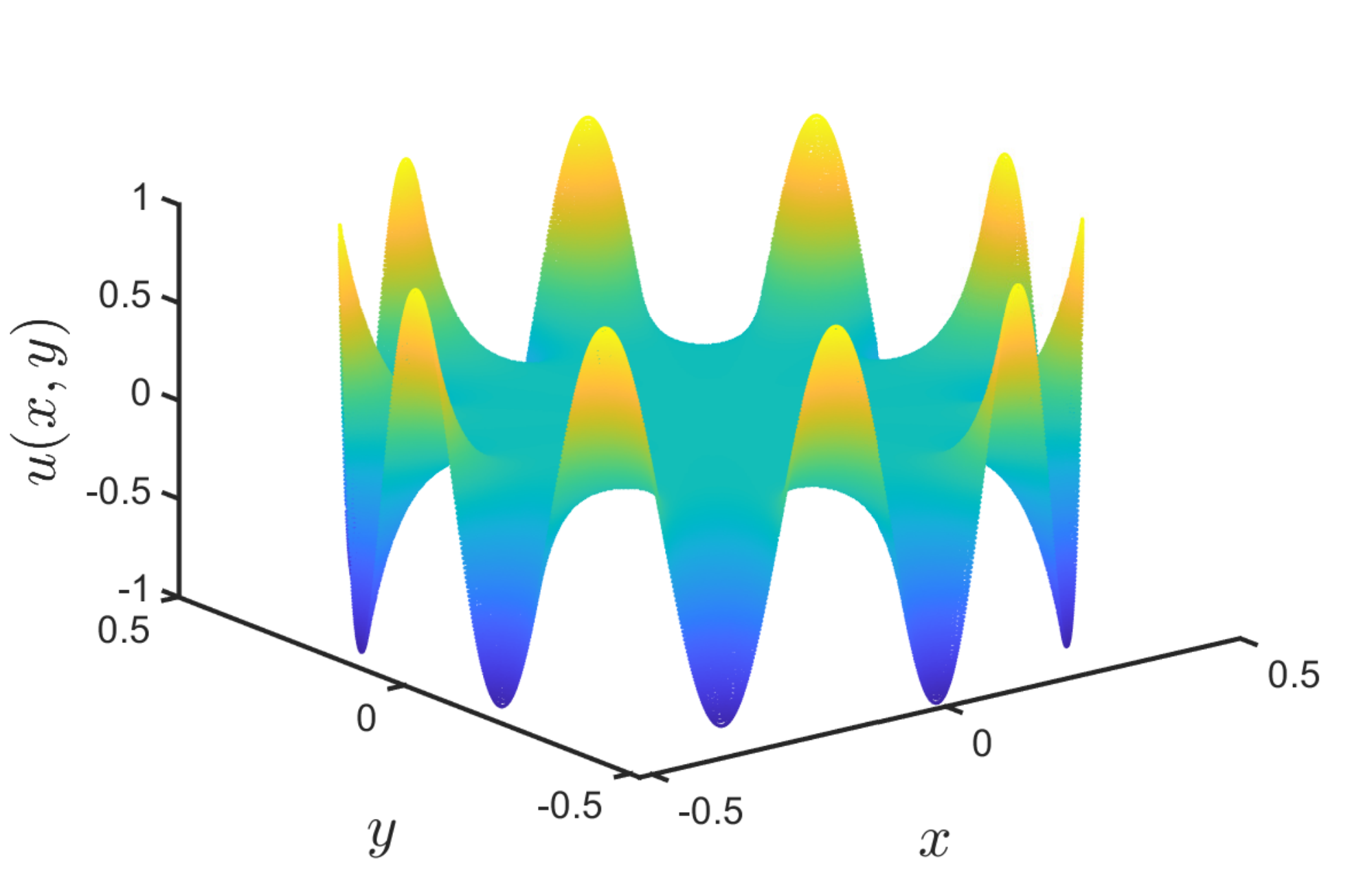}
        \caption{$u(r,\theta) = 1024\cos(10 \theta) r^{10}$}
    \end{subfigure} 
    \caption{The analytical solutions used for the (a) Poisson and (b) Laplace test problems.}
    \label{fig:ml0}
\end{figure}

Traditionally, multigrid solvers have been used to accelerate the solution of large systems of equations \cite{BriggsMultigrid,trottenberg2000multigrid}. Multigrid methods, however, require structured grids as the name indicates. In this section, it is illustrated how to utilize the geometric flexibility of RBF-FD in combination with the proposed node subsampling strategy to setup a geometric multilevel solver. Each example uses spatially varying node sets that seek to match the solutions. Consider the two-dimensional Poisson problem,

\begin{equation} \label{eq:peq}
\begin{aligned} 
\mathcal \nabla^2 u &= f,\hspace{0.25cm} \bm{x} \in \Omega\\
u &= g, \hspace{0.25cm} \bm{x} \in \partial\Omega 
\end{aligned}
\end{equation}

\noindent where $u = u(\bm{x}) = u(x,y) \in \mathbb{R}$ is the exact solution in a disk with unit diameter, i.e., $\Omega = \{(x,y), x^2 + y^2 \leq 0.5\}$, $g = g(\bm{x}) \in \mathbb{R}$ specifies Dirichlet boundary conditions on $\partial\Omega$ and $f = f(\bm{x}) \in \mathbb{R}$ specifies the source term. Two different problems are solved using variable density node sets in order to test the applicability of the proposed node subsampling strategy. The first problem considered is the Poisson problem for which $g=0$ and $f= 200 e^{-100r} \left(100r-1\right)/r$ such that the solution given in polar coordinates is $u(r,\theta) = 2\exp(-r/ 0.01)$, while the other is a Laplace problem (i.e. $f = 0$) for which $g=\cos(10\theta)$ such that the solution given is $u(r,\theta) = 1024\cos(10 \theta) r^{10}$ (see Figure \ref{fig:ml0}).

%\begin{figure}[H] 
%\centering
%\includegraphics[width=0.65\linewidth]{Figures/ml0.pdf} 
%\caption{The analytical solution $u(\bm{x}) = \exp(-||\bm{x}||_2^2 / 0.01)$ in $\Omega$.}
%\label{fig:ml0}
%\end{figure}

%\begin{figure}[H] 
%\centering
%\includegraphics[width=0.65\linewidth]{Figures/ml0.pdf} 
%\caption{The analytical solution $u(\bm{x}) = \exp(-||\bm{x}||_2^2 / 0.01)$ in $\Omega$.}
%\label{fig:ml0}
%\end{figure}

%%%%%%%%%%%%%%%%%%%%%%%%%%%%%
\subsection{Radial basis function-generated finite differences}

To discretize the problem in (\ref{eq:peq}) using RBF-FD \cite{bayona2017role,FLYER201639,FF2015RBFPrimer} we introduce the polyharmonic spline (PHS) radial basis functions $\phi_i(r)= \phi \left(\|\bm{x} - \bm{x}_i\|_2 \right) = \|\bm{x} - \bm{x}_i\|_2^{2k+1}$ and bivariate monomials $p_j(\bm{x})$ to approximate the exact solution to (\ref{eq:peq}) as

\begin{equation} \label{eq:ansatz}
    u(\bm{x}) \approx u_h(\bm{x}) =  \sum_{i=1}^{n} \kappa_i \phi \left(\|\bm{x} - \bm{x}_i\|_2 \right) + 
    \sum_{j=1}^{\ell} \gamma_j p_j(\bm{x})
\end{equation}

\noindent which we require to match at $n$ nodes, i.e., spatial points $\{\bm{x}_i\}_{i=1}^{n}$,

\begin{equation}
u(\bm{x}_i) = u_h(\bm{x}_i), \hspace{0.25cm} \text{for} \hspace{0.25cm} i=1, 2, ..., n,
\end{equation}

\noindent while enforcing the additional constraints,

\begin{equation}
\sum_{i=1}^{n} \kappa_i p_j(\bm{x}_i) = 0, \hspace{0.25cm} \text{for} \hspace{0.25cm} j=1, 2, ..., \ell, 
\end{equation}

\noindent where $\ell = (m+1)(m+2)/2$ is the number of monomial terms in a bivariate polynomial of degree $m$. The above equations can be arranged in a linear system of equations,

\begin{equation} 
\tilde{A}     \begin{bmatrix}
        \bm{\kappa} \\
        \bm{\gamma}
    \end{bmatrix}
    =
    \begin{bmatrix}
        A & P\\
        P^T & 0
    \end{bmatrix}
    \begin{bmatrix}
        \bm{\kappa} \\
        \bm{\gamma}
    \end{bmatrix}
    =
    \begin{bmatrix}
        \bm{u} \\
        \bm{0}
    \end{bmatrix},
\end{equation}

\noindent where $\bm{\kappa}, \bm{u} \in \mathbb{R}^n$, $\bm{\gamma} \in \mathbb{R}^\ell$, $A_{ij} = \phi \left(\|\bm{x}_i - \bm{x}_j\|_2 \right)$ is an entry of the RBF collocation matrix $A \in \mathbb{R}^{n \times n}$ and $P_{ij}= p_j(\bm{x}_i)$ is an entry of the supplementary polynomial matrix $P \in \mathbb{R}^{n\times \ell}$. Now, the linear operation $\mathcal{L}$ can be approximated at an evaluation point $\bm{x}_e$ as,

\begin{equation}
\mathcal{L} u|_{\bm{x}_e} \approx  \sum_{i=1}^{n} \kappa_i \mathcal{L}\phi \left(\|\bm{x} - \bm{x}_i\|_2 \right)|_{\bm{x}_e} + 
    \sum_{j=1}^{\ell} \gamma_j \mathcal{L} p_j(\bm{x})|_{\bm{x}_e},
\end{equation}

\noindent which again can be arranged in matrix-vector format,

\begin{equation}
   \mathcal{L} u|_{\bm{x}_e}  \approx  
   \left[  \bm{a}^T \hspace{0.2cm} \bm{b}^T \right]     
   \begin{bmatrix}
        \bm{\kappa} \\
        \bm{\gamma}
    \end{bmatrix}
    =
    \left[  \bm{a}^T \hspace{0.2cm} \bm{b}^T  \right]  \tilde{A}^{\text{-}1} 
        \begin{bmatrix}
        \bm{u} \\
        \bm{0}
    \end{bmatrix}
    =
       \left[  \bm{w}^T \hspace{0.2cm} \bm{v}^T \right]     
        \begin{bmatrix}
        \bm{u}\\
        \bm{0}
    \end{bmatrix}
    = \bm{w}^T  \bm{u}
\end{equation}

\noindent where $a_{i} = \mathcal{L}\phi \left(\|\bm{x} - \bm{x}_i\|_2 \right)|_{\bm{x}_e}$ corresponds to the $i$th entry of $\bm{a} \in \mathbb{R}^{n}$ and $b_{j} = \mathcal{L} p_j(\bm{x})|_{\bm{x}_e}$ corresponds to the $j$th entry of $\bm{b} \in \mathbb{R}^{\ell}$. The weights necessary for the multilevel solver, i.e. $\bm{w} \in \mathbb{R}^{n}$, can equivalently be computed by solving the linear system,

\begin{equation} 
    \begin{bmatrix}
        A & P\\
        P^T & 0
    \end{bmatrix}
    \begin{bmatrix}
        \bm{w} \\
        \bm{v}
    \end{bmatrix}
    =
    \begin{bmatrix}
        \bm{a} \\
        \bm{b}
    \end{bmatrix}.
\end{equation}

%%%%%%%%%%%%%%%%%%%%%%%%%%%%%
\subsection{Geometric multilevel elliptic solver}

The meshfree geometric multilevel solver introduced here is based on similar ideas as the meshfree geometric multilevel method \cite{WrightJonesShankarMGMRBFFD} used for solving PDEs on surfaces, although some parts differ. In this study, no Krylov subspace methods will be used to increase the rate of convergence \cite{krylovghai}. Furthermore, the coarse grid difference operators are computed explicitly on each node set level. Finally, all restriction operations are performed as injection (i.e. directly using values from the fine node set). For further details on multilevel methods, the reader is referred to literature on the topics of multilevel approximation \cite{fasshauer2007meshfree,wendland2004scattered} and multilevel solvers \cite{WrightJonesShankarMGMRBFFD,zamolo2019novel}. 

A pseudocode for the proposed geometric multilevel solver is given in Algorithm \ref{alg:mlsolver}, where  $u_1$ is the solution at the finest node set level, $L = \{L_j\}_{j=1}^{p}$ are the difference operators computed for each level, $I = \{I_{j+1}^{j}\}_{j=1}^{p-1}$ are the interpolation operators for each level and $R = \{R_{j}^{j+1}\}_{j=1}^{p-1}$ are the restriction (injection) operators for each level. The multilevel solver performs up to $i_{max}$ iterations (V-cycles) unless the relative residual $||r_1^i||_2 / ||r_1^0||_2 = ||f_1 - L_1 u_1^i||_2 / ||f_1 - L_1 u_1^0||_2$ of the $i$th iteration becomes less than a predefined tolerance $tol$.

The basis of the geometric multilevel solver is the geometric multilevel V-cycle, which is described in Algorithm \ref{alg:mlvcyc}. During the V-cycles, pre- and post smoothing operations are performed using $(\nu_1,\nu_2)$ Gauss-Seidel relaxations, respectively. At the coarsest node set level, a sparse LU solver is used.

The pseudocode for performing the geometric multilevel preprocessing, i.e., establishing all the different subsets of nodes, $X = \{X_j\}_{j=1}^p$, and the discrete operators $L, I, R$, is given in Algorithm \ref{alg:mlpre}. The necessary input for this algorithm are two node sets, $\{X_{bg},X_{b} \}$, which describe the scattered node set covering $\Omega$ and boundary nodes on $\partial\Omega$ at the finest node set level, respectively. Finally, the parameter $N_{\text{min}}$ is used to control the minimum number of boundary nodes at the coarsest node set level.

\begin{algorithm}
  \caption{Geometric multilevel solver} \label{alg:mlsolver}
  \begin{algorithmic}[1]
  \Function{mlsolver}{$u_1, f_1, L, I, R, \nu_1, \nu_2, tol, i_{max}$}
    \While{$i < i_{max}$ and $||r_1||_2 < ||f_1||_2 \cdot tol$}
      \State $i \leftarrow i + 1$
        \State $u_1 \leftarrow \textproc{mlvcyc} (u_1, f_1, L, I, R, \nu_1, \nu_2)$
    \EndWhile
    \State \Return u
  \EndFunction
  \end{algorithmic}
\end{algorithm}

\begin{algorithm}
  \caption{Geometric multilevel V-cycle} \label{alg:mlvcyc}
  \begin{algorithmic}[1]
  \Function{mlvcyc}{$u_1, f_1, L, I, R, \nu_1, \nu_2$}
    \State $u_1 \leftarrow \textproc{relax}(u_1,f_1,L_1,\nu_1)$
    \State $r_2 \leftarrow R_1^2 (f_1 - L_1 u_1)$
    \For{$j = 2$ to $p-1$}
    \State $e_j \leftarrow \textproc{relax}(0,r_j,L_j,\nu_1)$
    \State $r_{j+1} \leftarrow R_j^{j+1} (r_j - L_j e_j)$
    \EndFor
    % \State $e_p = L_p^{\text{-}1} r_p$
    \State $e_p = \textproc{lusolve} (L_p,r_p)$
    \For{$j = p-1$ to $2$}
    \State $e_{j} \leftarrow e_{j} + I_{j+1}^j e_{j+1}$
    \State $e_j \leftarrow \textproc{relax}(e_j,r_j,L_j,\nu_2)$
    \EndFor
    \State $u_1 \leftarrow u_1 + I_{2}^1 e_{2}$
    \State $u_1 \leftarrow \textproc{relax}(u_1,f_1,L_1,\nu_2)$
    \State \Return $u_1$
  \EndFunction
  \end{algorithmic}
\end{algorithm}

\begin{algorithm}
  \caption{Geometric multilevel preprocessing} \label{alg:mlpre}
  \begin{algorithmic}[1]
    \Function{mlpre}{$X_{bg}, X_b, N_{\text{min}}$}
    \State $[X, R, p] \leftarrow \textproc{mlmfsub} (X_{bg},X_{b},N_{\text{min}})$
            \State $L_1 \leftarrow \textproc{rbffd} (X_1,X_1)$
     \For{$j=1$ to $p-1$}
        \State $I_{j+1}^{j} \leftarrow \textproc{rbffd} (X_{j+1},X_{j})$
        \State $L_{j+1} \leftarrow \textproc{rbffd} (X_{j+1},X_{j+1})$
     \EndFor
     \State \Return $X, L, I, R, p$
  \EndFunction
  \end{algorithmic}
\end{algorithm}

% \begin{algorithm}
%   \caption{Moving front node subsampling} \label{alg:mlnsubs}
%   \begin{algorithmic}[1]
%   \Function{nsubs}{$u_1, f_1, L, I, R, \nu_1, \nu_2, tol, i_{max}$}
%   \State $j = 1$
%   \State $f_c = 1.5$
%   \State $K = 10$
%     \While{size$(X_j,1) \geq N_{min}$}
%         \State $X_j \leftarrow \textproc{nsubs} (X_j, f_c, K)$
%         \State $X_j \leftarrow \textproc{nsubs} (X_j, f_c, K)$

%     \State $j \leftarrow j + 1$

%     \EndWhile
%     \Return u
%   \end{algorithmic}
% \end{algorithm}

%%%%%%%%%%%%%%%%%%%%%%%%%%%%%
% \subsection{Test problems}\label{sec:ExampleNumRes}

\begin{figure}[H] 
\includegraphics[width=\linewidth]{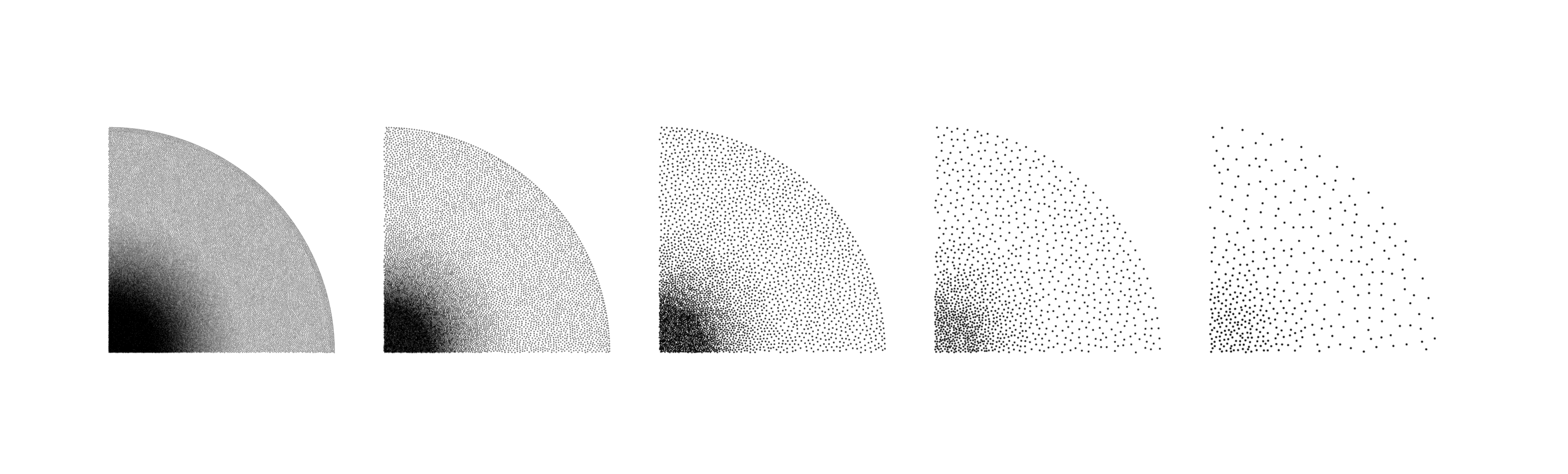} 

\vspace{-1.25cm}

\includegraphics[width=\linewidth]{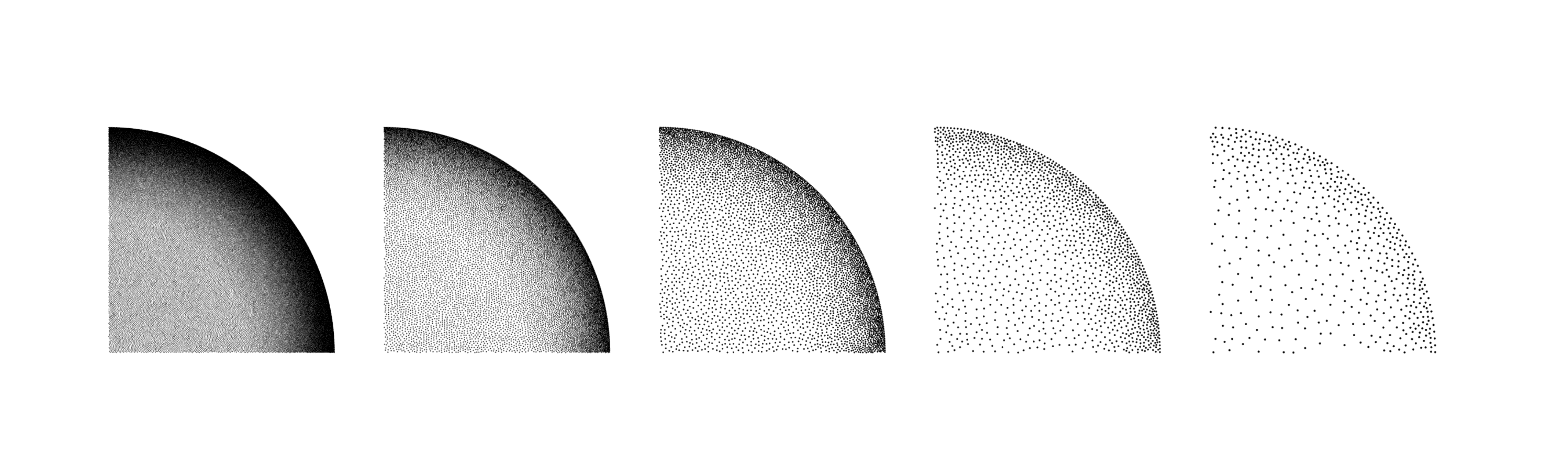} 
\vspace{-1.25cm}
\caption{Example of the multilevel node subsampling process for the Poisson problem node set (top) and the Laplace problem node set (bottom). Only nodes within the first quadrant of the Cartesian coordinate system are shown.}
\label{fig:ml1}
\end{figure}

The multilevel solver is tested on node sets that have been generated with variable node densities as illustrated in Figure \ref{fig:ml1}, where $N_{\text{min}} = 60$ for the Poisson problem and $N_{\text{min}} = 120$ for the Laplace problem. The node density function used in this study is defined by a linear transition between two prescribed node densities as,

\begin{equation}
\rho(d)
= 
\begin{cases}
    \rho_{\text{1}},& d < d_{\text{lim}}\\
    \rho_{\text{1}} + (\rho_{\text{2}}-\rho_{\text{1}}) (d-d_{\text{lim}})/d_{\text{bl}},& d_{\text{lim}} \leq d \leq d_{\text{lim}} + d_{\text{bl}}\\
    \rho_{\text{2}},              & \text{otherwise}
\end{cases}
\end{equation}

\noindent where $d = || \bm{x} ||_2$ is the distance to the origin according to Figure \ref{fig:ml0}, $\rho_{\text{1}}$ is the node density in region 1, $d_{\text{lim}}$ is a distance within which $\rho_{\text{1}}$ is kept constant, whereas $d_{\text{bl}}$ is the distance over which $\rho_{\text{1}}$ linearly blends into $\rho_{\text{2}}$. It should be noted that the nodes in the vicinity of boundary have been adjusted by means of repulsion only at the finest node set level \cite{FF2015NodeGen}. It should be noted that the node density function, $\rho(d)$, is chosen such that node density can be matched with the characteristics of the solutions.

% \begin{figure}[H] 
% \includegraphics[width=\linewidth]{Figures/ml1.pdf} 
% \caption{Example of the multilevel node subsampling process for $N = 56870$, which provides subsampled node sets equal to $N = \{56870, 16453, 5223\}$. Here the variable node density is generated based on the node density function $\rho(\bm{x}) = \rho_b/(1 + 2u(\bm{x}))$, where $\rho_b$ is the node density on $\partial \Omega$.}
% \label{fig:ml1}
% \end{figure}

The numerical test setups have been chosen to showcase the applicability of the node subsampling strategy for multilevel solvers using RBF-FD and to test whether the high-order accuracy will still be dictated by the degree of the augmented polynomials as shown, e.g., in \cite{bayona2017role}. Thus, the parameters used for computing the difference operators, $L$, of polynomial degree $m_{L}$ are chosen to be $(k,n) = (1,2\ell)$, while the parameters $(m_{I},k,n) = (0,0,5)$ are used for computing the interpolation operators, $I$. The parameters for the interpolation operators are kept fixed for all choices of $m_{L}$. Finally, the multilevel solver settings are defined as $(\nu_1,\nu_2,i_{max},tol) = (2,1,50,10^{-16})$ for both test problems, whereas the polynomial degree $m_L$ ranges from 2 to 8.

\subsection{Poisson Equation Test Problem}\label{sec:ExampleNumRes}

\begin{figure}[H] 
\includegraphics[width=\linewidth]{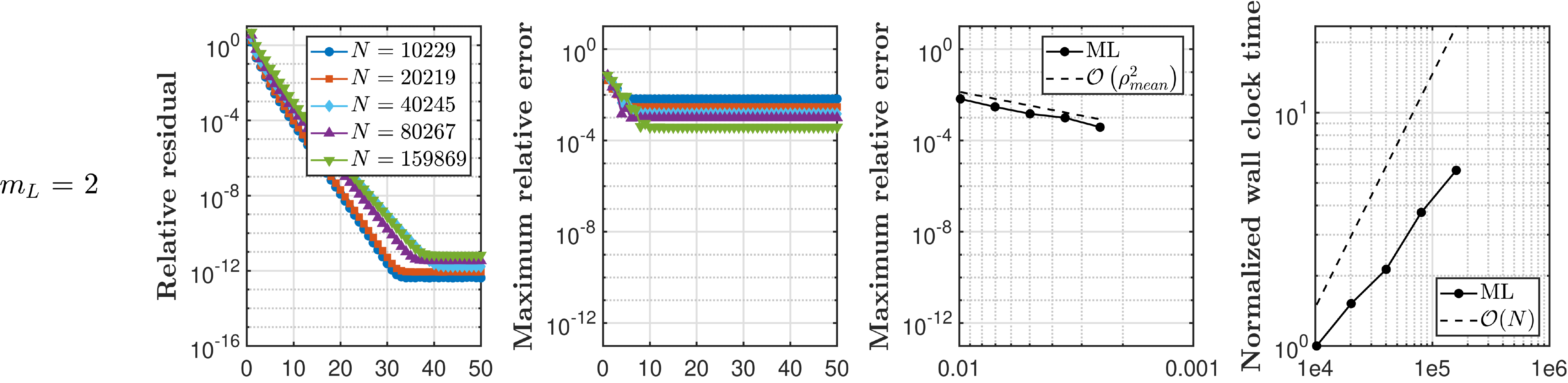} \\

\includegraphics[width=\linewidth]{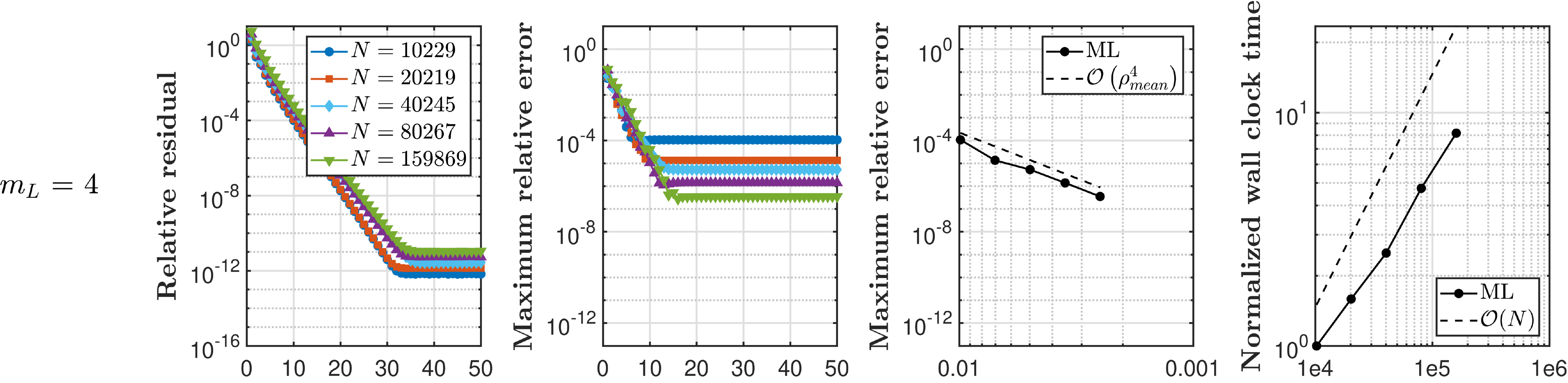}\\

\includegraphics[width=\linewidth]{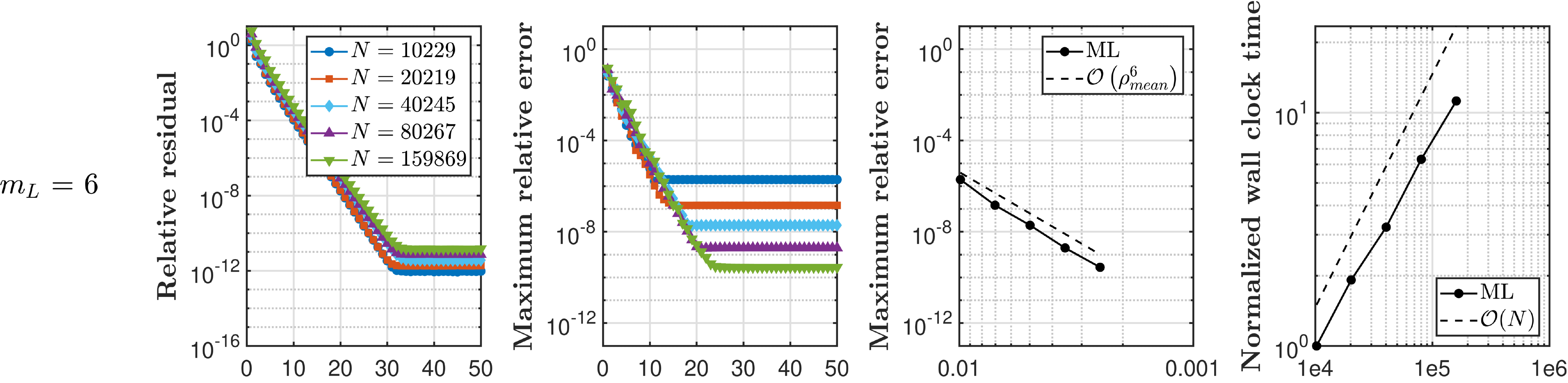}\\

\includegraphics[width=\linewidth]{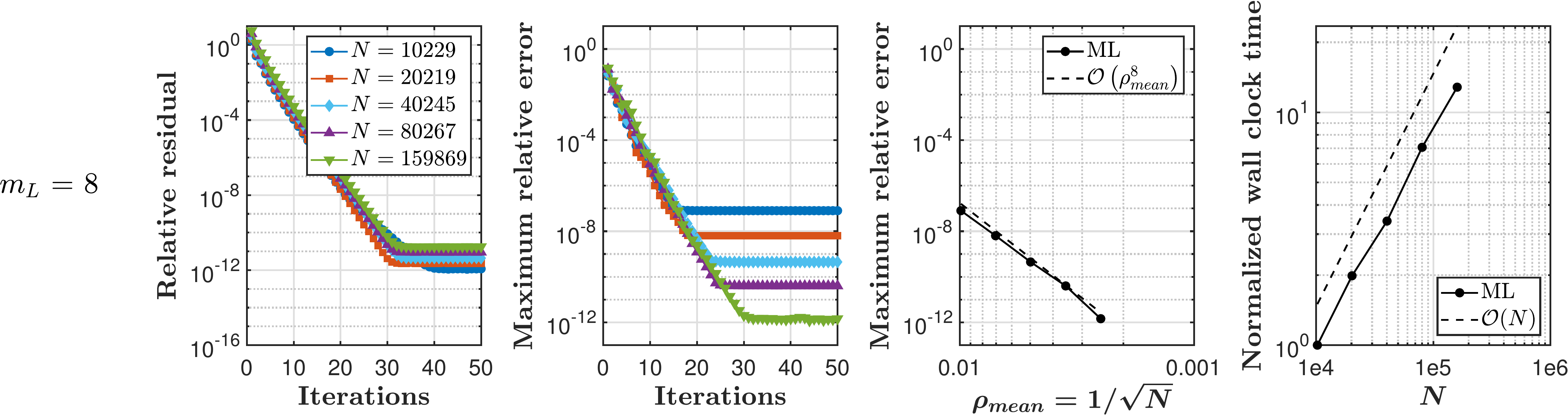}
\caption{Poisson problem performance indicators of the implemented geometric multilevel solver for polynomial degrees of $m_L = \{2,4,6,8\}$ from top to bottom. The mean node density is defined as $\rho_{mean} = 1/\sqrt{N}$.}
\label{fig:ml2_poisson}
\end{figure}

First, it can be seen from Figure \ref{fig:ml2_poisson} that the wall clock time scales linearly with the number of nodes, which is in accordance with expectations for any multigrid or multilevel solver \cite{BriggsMultigrid,trottenberg2000multigrid,WrightJonesShankarMGMRBFFD}. Furthermore, the maximum relative error ($||u-u_h||_{\infty}/||u||_{\infty}$) decreases as function of node set resolution ($\rho_{mean} = 1/\sqrt{N}$) and the slope is dictated by the polynomial degree of the difference operator, $m_L$. This is in agreement with previous RBF-FD studies \cite{bayona2017role}.

In this study, the low-order interpolation operators are chosen in order to accelerate the convergence of the multilevel solver. However, this choice is not aligned with the rule of thumb for the transfer operators, i.e. $m_I + m_R > 2$, which is used in traditional multigrid methods for the Poisson problem \cite{trottenberg2000multigrid}. In this work, the order of the restriction operators are $m_R = 0$ because injection is used as the restriction operation between all node set levels. Nevertheless, no detrimental effects have been noticed in any of the numerical tests conducted.

\begin{figure}[H]
\includegraphics[width=\linewidth]{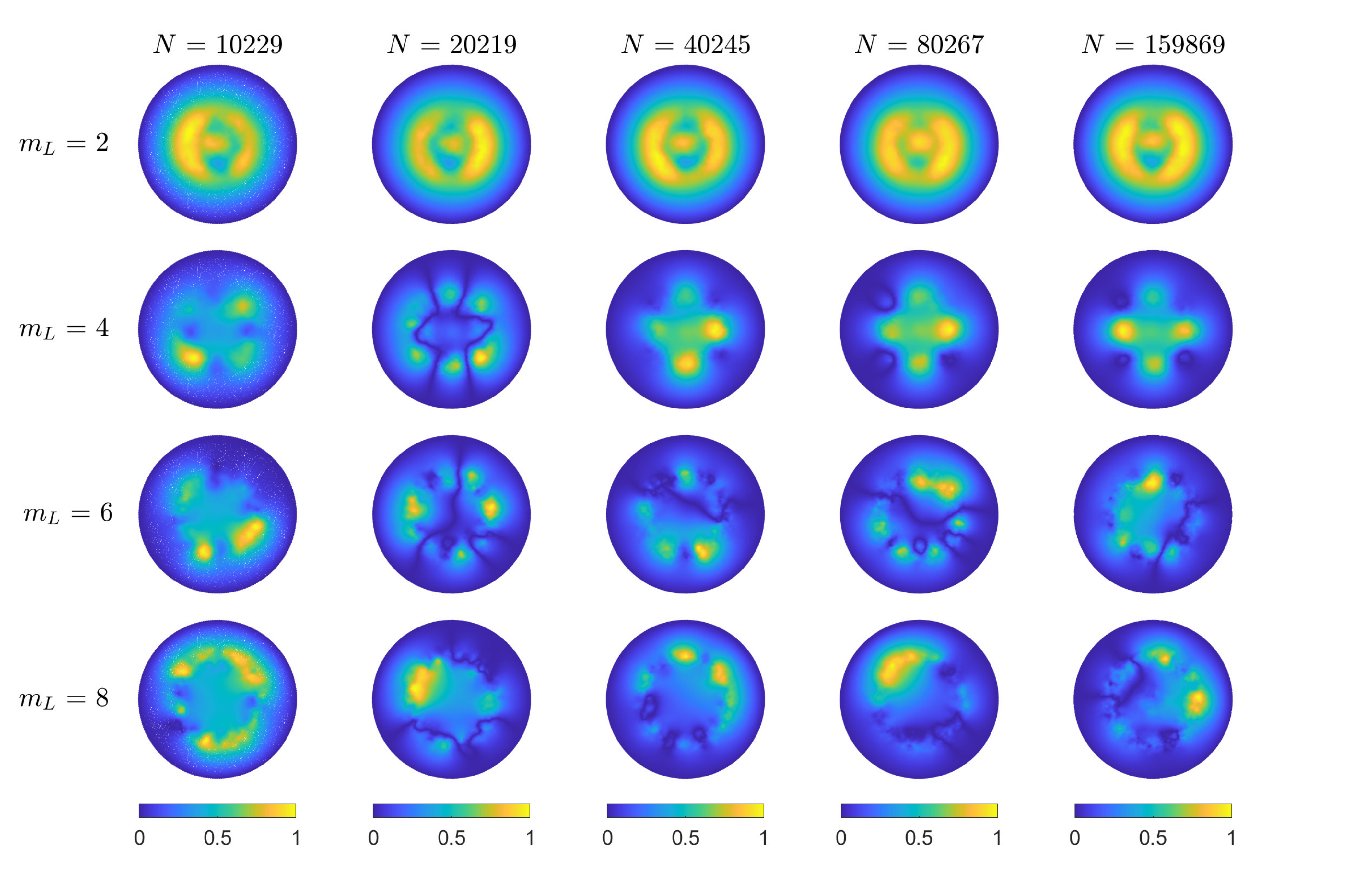}
\caption{Normalized relative error distributions for various orders of the difference operators and node set resolutions for the Poisson problem. The color scale factors refer to the plateau of the maximum relative error plots in figure \ref{fig:ml2_poisson}.}
\label{fig:ml3_poisson}
\end{figure}

Finally, the proposed multilevel solver should provide solutions without any directional bias. Hence, to identify whether any directional bias is present, the relative errors for all node set resolutions have been normalized and depicted in Figure \ref{fig:ml3_poisson}. Thus, the scale factors used in Figure \ref{fig:ml3_poisson} refer to the plateau of the maximum relative error plots in Figure \ref{fig:ml2_poisson}. As no particular directional pattern is noticed in Figure \ref{fig:ml3_poisson}, it can be concluded that the subsampling process used for setting up the multilevel solver does not introduce any  directional bias.

%%%%%%%%%%%%%%%%%%%%%%%%%%%%%
\subsection{Laplace Equation Test Problem}

\begin{figure}[H] 
\includegraphics[width=\linewidth]{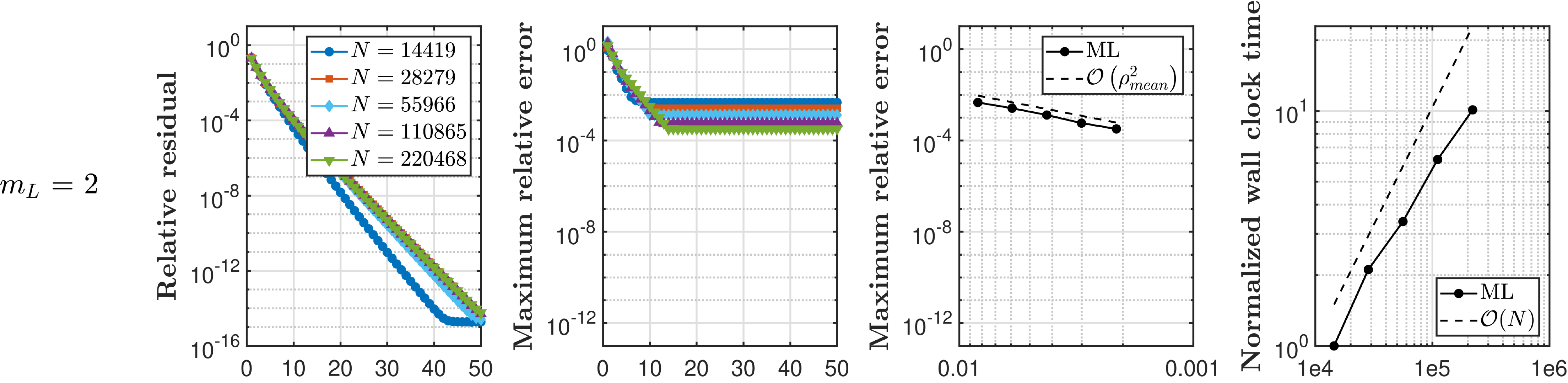} \\

\includegraphics[width=\linewidth]{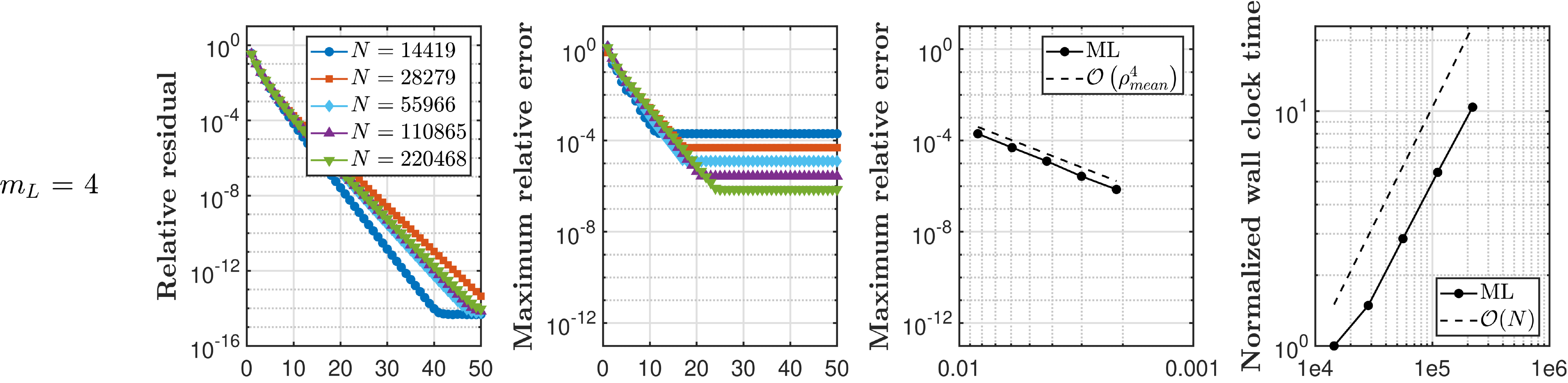}\\

\includegraphics[width=\linewidth]{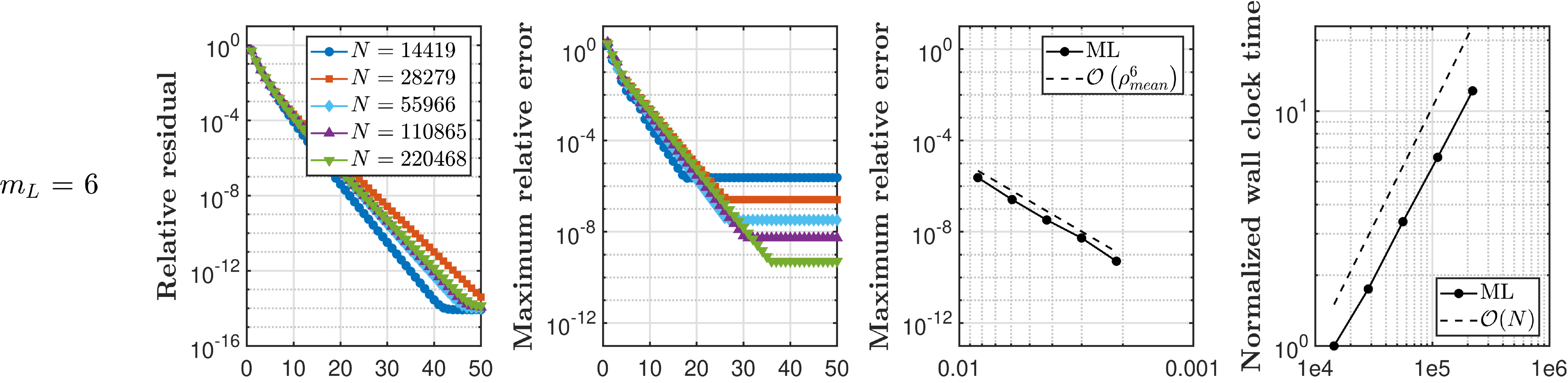}\\

\includegraphics[width=\linewidth]{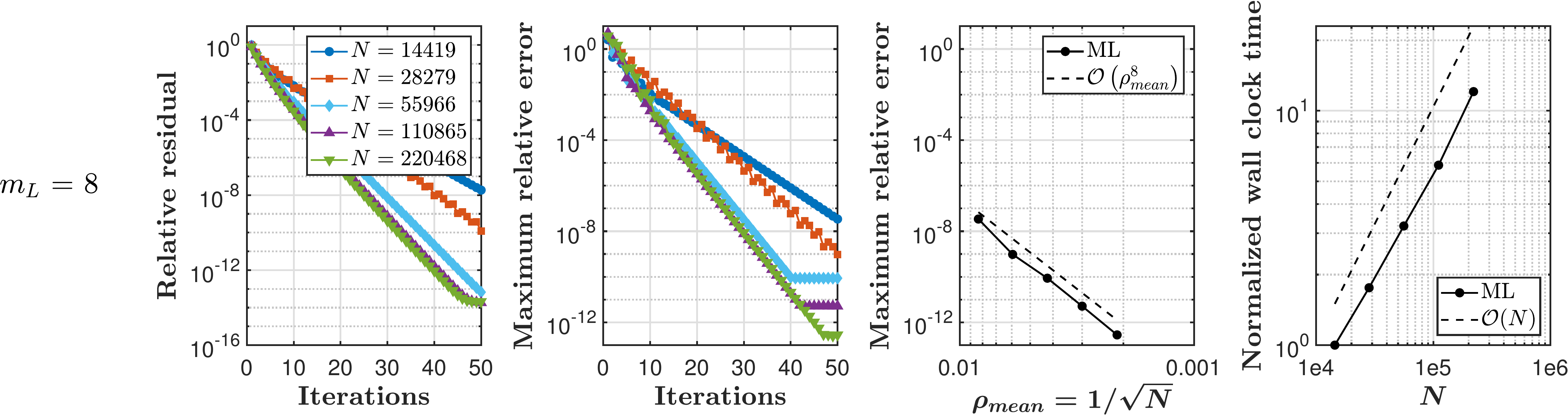}
\caption{Laplace problem performance indicators of the implemented geometric multilevel solver for polynomial degrees of $m_L = \{2,4,6,8\}$ from top to bottom. The mean node density is defined as $\rho_{mean} = 1/\sqrt{N}$.}
\label{fig:ml2_laplace}
\end{figure}

\begin{figure}[H]
\includegraphics[width=\linewidth]{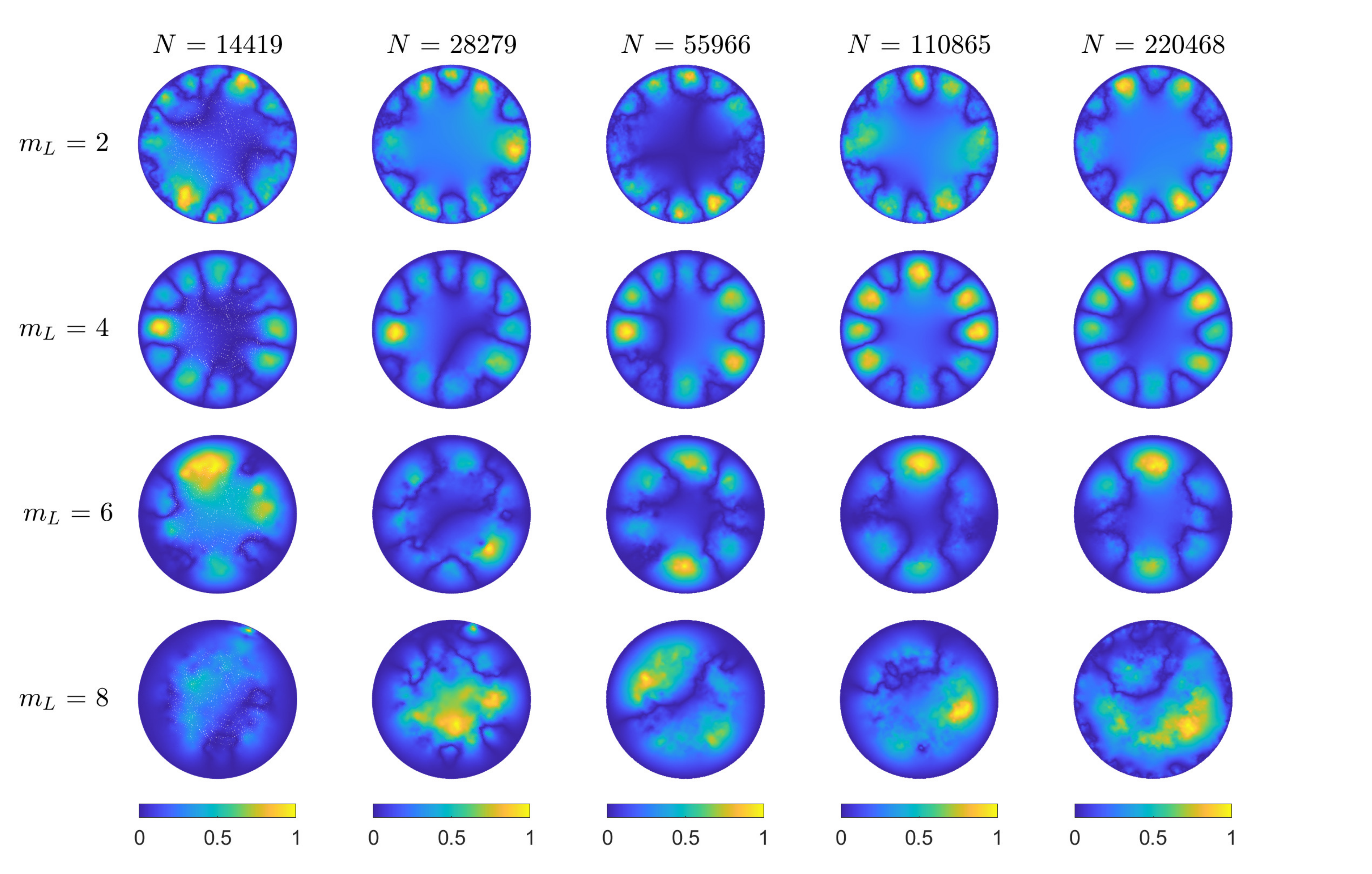}
\caption{Normalized relative error distributions for various orders of the difference operators and node set resolutions for the Laplace problem. The color scale factors refer to the plateau of the maximum relative error plots in figure \ref{fig:ml2_laplace}.}
\label{fig:ml3_laplace}
\end{figure}

The performance indicators of the multilevel solver for the Laplace problem are illustrated in Figure \ref{fig:ml2_laplace}. The same overall conclusions that were made for the Poisson problem can be made for the Laplace problem, i.e. high-order accuracy and linear scaling of the computation time. For $m_L = 8$, note that the convergence factor ($||r_1^i||_2/||r_1^{i-1}||_2$) of the multilevel solver is less for $N = 14419$ and $N = 28279$ compared with the other values of $N$. The decrease in convergence factors is most likely caused by a stencil size that is too large as compared to the relatively low node density near the boundary since localized error peaks are present for both $N = 14419$ and $N = 28279$ ($m_L = 8$) in Figure \ref{fig:ml3_laplace}. This decrease in convergence factor does not occur if the node resolution is increased to $N = 55966$ or above. Furthermore, if the multilevel V-cycle is used as a preconditioner for a Krylov subspace method, e.g. the generalized minimal residual method or biconjugate gradient stabilized method, fewer iterations will be needed for the solution to converge and the solver will be more robust compared to the standalone multilevel solver. However, each iteration will become more computational expensive. Thus, whether the multilevel V-cycle should be used as a standalone solver or a preconditioner is a trade-off between computational cost and robustness.

The normalized relative error distributions in Figure \ref{fig:ml3_laplace} illustrate that no directional bias seems to be introduced by the subsampling process, which is the same conclusion as for the Poisson problem.

%%%%%%%%%%%%%%%%%%%%%%%%%%%%%%%%%%%%%%%%%%%%%%%%%%%%%%%%%%
\section{Conclusion}
A novel method for subsampling quasi-uniform node sets of highly variable density with sharp gradients is presented along with boundary preservation techniques and two novel measures for evaluating node quality of subsamplings. The moving front subsampling algorithm demonstrates the capability to coarsen a node set with high contrast and detail. Additionally, the moving front algorithm maintains the characteristics of the original node set as outlined by the comparative local regularity of the average distance and standard deviation of distances to the $k$ nearest nodes. It is also faster, both by a constant and in the limit as node set size increases, than any other algorithm considered in this paper for subsampling variable density node sets. 

The utility of the moving front algorithm for the purpose of subsampling node sets in a meshfree multilevel PDE solver is demonstrated by solving both the Poisson and Laplace problems on variable density node sets. In both test cases, the meshfree PDE solver with the multilevel method and the proposed subsampling algorithm achieves the fast linear scaling of computational cost with node set size expected from a multilevel scheme. At the same time, this combination has no adverse impact on the expected high-order accuracy of the RBF-FD method. The meshfree multilevel PDE solver has been tested up through eighth order convergence and also demonstrates very robust performance.\\

\acknowledgments{Andrew Lawrence acknowledges support from the US Air Force\footnote{The views expressed in this article are those of the authors and do not reflect the official policy or position of the
Air Force, the Department of Defense or the U.S. Government.}}

\pagebreak
\appendix
%%%%%%%%%%%%%%%%%%%%%%%%%%%%%%%%%%%%%%%%%%%%%%%%%%%%%%%%%%
\section{Subsampling Algorithms}

%%%%%%%%%%%%%%%%%%%%%%%%%%%%%
\subsection{Moving Front Subsampling}\label{apx:MF}

The Python code for the moving front subsampling algorithm. The code can also be found on the author's GitHub repository in both MATLAB and Python along with examples of implementation \cite{Lawrence_GitHub}.

\begin{verbatim}

import numpy as np
from sklearn.neighbors import NearestNeighbors

def MFNUS(xy, fc=1.5, K=10):
    """
    Moving Front Non-Uniform Subsampling
    
    Args:
        xy (array): initial node set to be subsample
        c (float): coarsening factor
        K (float): number of nerarest neighbors to check in algorithm
        
    Returns:
        xy_sub (array): subsampled node set
    """
    if xy.shape[0] < xy.shape[1]:
        xy = xy.T

    # algorithm
    N = xy.shape[0]  # Get the number of its dots
    sort_ind = np.lexsort(xy.T,axis=0)
    xy = xy[sort_ind, :]  # Sort dots from bottom and up
    
    # Create nearest neighbor pointers and distances
    nbrs = NearestNeighbors(n_neighbors=K+1, algorithm='auto').fit(xy)
    distances, indices = nbrs.kneighbors(xy)
    
    for k in range(N):  # Loop over nodes from bottom and up
        if indices[k, 0] != N+1:  # Check if node already eliminated
            ind = np.where(distances[k, 1:] < fc*distances[k, 1])[0]
            ind2 = indices[k, ind+1]
            ind2 = np.delete(ind2,ind2 < k)   # Mark nodes above present one, and which
            indices[ind2, 0] = N+1        # are within the factor fc of the closest one

    elim_ind_sorted = indices[:, 0] != N+1
    xy_sub = xy[elim_ind_sorted]
        
    return xy_sub

\end{verbatim}

\bibliography{subsampling}
%\printbibliography

\end{document}